\documentclass{article}
\usepackage[top=2.5cm,bottom=2.5cm,left=2.0cm,right=2.5cm]{geometry}
\usepackage{tikz}
\usepackage{authblk}
\usepackage{amsmath}
\usepackage{amssymb}
\usepackage{amsthm}
\usepackage{color}
\usepackage[normalem]{ulem}
\usepackage{bm}
\usepackage{cases}
\usepackage[colorlinks,linkcolor=blue,anchorcolor=green,citecolor=red]{hyperref}
\usepackage{appendix}
\usepackage{graphicx}
\usepackage{subfigure}
\usepackage[section]{placeins}
\usepackage{booktabs}
\usepackage{lscape}
\usepackage[figuresright]{rotating}
\usepackage{cite}
\usepackage{algorithm}
\usepackage{algpseudocode} 
\usepackage{tikz,float}
\usepackage{caption}
\usepackage{bbding}
\DeclareUnicodeCharacter{00B1}{\ensuremath{\pm}}
\usetikzlibrary{arrows.meta}
\usetikzlibrary{patterns}
\allowdisplaybreaks[4]
\numberwithin{equation}{section}

\newtheorem{lemma}{Lemma}[section]
\newtheorem{theorem}{Theorem}[section] 

\newtheorem{remark}{Remark}[section] 

\title{Two-phase visco-elastic FENE-P model: Thermodynamically consistent modeling, sharp interface limit analysis and unconditionally energy stable numerical scheme}
\author[a]{Yuzhe Qin}
\author[b,*]{Dongdong He}
\affil[a]{School of Mathematics and Statistics,  
Shanxi University, Taiyuan 030006, Shanxi, China}
\affil[b]{School of Science and Engineering, 
The Chinese University of Hong Kong (Shenzhen), Shenzhen 518172, Guangdong, China}
\affil[*]{Corresponding author, hedongdong@cuhk.edu.cn}

\DeclareMathOperator{\tr}{tr}

\newcommand{\normh}[1]{\left\|#1\right\|_h}
\newcommand{\Rey}{\mathrm{Re}}
\newcommand{\Ca}{\mathrm{Ca}}
\newcommand{\Wi}{\mathrm{Wi}}
\newcommand{\Pe}{\mathrm{Pe}}
 
\begin{document}
\maketitle
\begin{abstract}

Although mathematical modeling and numerical simulations for two-phase Newtonian flows have been extensively studied for many years, relatively little work has focused on developing thermodynamically consistent models for two-phase Non-Newtonian viscoelastic flows, and to the author's best konowledge, no structure-preserving numerical methods have been reported for such models.  
In this paper, we focus on modeling the two-phase non-Newtonnian visco-elastic FENE-P model based on the energetic variational approach (EnVarA). 
In the model, a phase-field function is introduced to describe the diffuse interface which satisfies the Cahn-Hilliard (CH) equation. And all the  physical variables can be connected smoothly. We conduct a matched asymptotic analysis for the sharp interface limit, 
which shows the diffuse interface model will converge to the sharp interface model 
when the interface thickness trends to zero. 
We then reformulate the model in an equivalent form that facilitates the construction of the numerical scheme. A weakly decoupled scheme is proposed for the resulting FENE-P–CH system, and its unconditional energy stability is proved theoretically. Finally, a series of numerical experiments is conducted to validate the model and assess the performance of the proposed scheme.
\end{abstract}
\section{Introduction}
The two-phase and multiphase immiscible flows appear in many real-life science and engineering applications. These flows are central topics related to droplet formation, which is the natural phenomenon that appears in many industrial applications, such as ink jet printing~\cite{yasui2011progress}, coating~\cite{champougny2021dip,Abels2012Thermodynamically}, liquid spraying or atomization~\cite{nasr2013industrial}, and microfluidic devices~\cite{yeo2011microfluidic}.   It is essentially important to control the shape and size of the droplets in these applications, which are closely related to the quality of these industrial products. Polymeric fluids are widely used in these applications. Polymeric fluids usually exhibit complex non-Newtonian viscoelastic rheological properties which give rise to both elastic and viscous behaviors.

Mathematical modeling for two-phase and multiphase flows are helpful in understanding the underlying physics of aforementioned phenomenon. Although theoretical models for each phase fluid dynamics  are well established, mathematical complexity of these interfacial flow models greatly restricts analytical studies. Moreover, suitably modeling the motion of the fluid interface is necessary for two-phase and multiphase flows. All these factors give rise to the complexities and  call for numerical approaches.  

When numerical simulating two-phase and multiphase flow problems, in general, there are two types of methods for describing the moving interface, one type is the sharp interface method which considers the interface as zero thickness and imposes interface conditions directly, such as immersed boundary method~\cite{peskin2002immersed}, level set method~\cite{osher2001level}, immersed interface method~\cite{lee2003immersed}, volume of fluid~\cite{hirt1981volume}.  The other type is the diffusive interface method~\cite{anderson1998diffuse} which considers the interface as non-zero thickness, this method usually adopts a phase field variable defined globally in the whole domain and this phase field variable makes the smooth transition from one value to another. And phase field method~\cite{li2026global,huang2024maximum} usually adopts the variational approach to derive a thermodynamical consistent system, and the system satisfies the total energy decay property, and specially designed numerical methods associated with phase field method can preserve the energy dissipative law~\cite{el2022thermodynamically,feng2023energy}.  

For numerical simulation of two-phase Newtonian flows, there are numerous published work in the literature. Here, just name a few, Huang and Guo~\cite{huang2018improved} proposed an improved penalty immersed boundary method for multiphase flow simulation.  Lyras et al.~\cite{lyras2023conservative} used a conservative level set method for simulating liquid-gas flows and study the application in liquid jet atomisation. Thuy-Petrov et al.~\cite{thuy2025volume} developed the volume of fluid simulations of multiphase flows with high surface tension and curvature using a tensile force method.  Jain~\cite{jain2022accurate},  Chiu and Lin~\cite{chiu2011conservative} both proposed conservative phase-field methods for simulation of two-phase flows.  The fluids in these work are all Newtonian fluids.

Comparing to Newtonian fluid models, viscoelastic fluid models contain both viscous and elastic stresses. Different constitutive models have been proposed for the viscoelastic fluids, such as Oldroyd-B~\cite{oldroyd1950formulation}, Giesekus~\cite{giesekus1982simple}, Phan-Thien-Tanner (PTT)~\cite{thien1977new}, eXtended Pom-Pom (XPP)~\cite{verbeeten2001differential}, FENE-P~\cite{purnode1998polymer}. And numerical simulations for single-phase and two-phase non-Newtonian viscoelastic flow also received a lot of attention recently. For example, Alves et al.~\cite{alves2021numerical} provided a review for the numerical simulation of viscoelastic fluid flows with an emphasis on the finite-volume method. López-Herrera et al.~\cite{lopez2019adaptive}  proposed an adaptive numerical solver for viscoelastic 2D two-phase flows based on the volume-of-fluid method, where the algorithms have been developed for  Oldroyd-B, FENE-P and PTT viscoelastic fluid models. Izbassarov and Muradoglu~\cite{izbassarov2015front} developed a front-tracking method for viscoelastic two-phase systems  in which one or both phases could be viscoelastic. Yu et al.~\cite{yu2007two} developed the two-phase level set method for simulating the Oldroyd-B viscoelastic ink jet pinch-off phenomena. Zografos et al.~\cite{zografos2020viscoelastic} developed a two-phase viscoelastic fluid solver based on the phase field approach. Doherty et al.~\cite{doherty2024log} developed numerical solvers for both single- and multi-phase axisymmetric viscoelastic flows.  Mokbel et al.~\cite{mokbel2018phase} derived the two-phase viscoelastic flow model and performed numerical simulations for large deformation of droplets.

It is well known that numerically preserving the original physical properties when numerically solving differential equations arising from different physical problems are greatly important. For instance, one active research direction is to design structure-preserving schemes for ordinary differential equations (see~\cite{lubich2006geometric} and references therein). 
Li and Vu-Quoc~\cite{li1995finite} pointed out that ``in some areas, the ability to preserve some invariant properties of the original differential equation is a criterion to judge the success of a numerical simulation". 
Many important structure preserving numerical methods have been proposed for various  partial differential equation models. 
For example, Shen and Yang~\cite{shen2010phase} constructed energy stable schemes for solving two-phase Allen--Cahn Navier--Stokes phase model with variable densities. 
Based on convex splitting technique, Chen et al.~\cite{chen2017uniquely} designed decoupled energy stable schems for the Cahn--Hilliard Stokes--Darcy system for two phase
flows in karstic geometry. 
Yang and Yu~\cite{yang2018efficient} developed second order unconditionally energy stable schemes for the phase field moving contact line model by using invariant energy quadratization method. 
Li and Shen~\cite{li2020sav} constructed the second-order energy stable method for Cahn-Hilliard Navier-Stokes phase-field model based on the scalar auxiliary variable approach. 
Hou et al.~\cite{hou2024energy} developed energy stable spectral renormalization exponential integrator method for gradient flow problems. Li et al.~\cite{li2026second} developed the unit length preserving and energy dissipative Runge-Kutta method for gradient flows of vector fields. Duan et al.~\cite{duan2026third} provided structure-preserving exponential time differencing Runge--Kutta scheme for the fluid-surfactant phase field Model. 
And structure preserving methods have also been developed for Poisson--Nernst--Planck and Navier--Stokes--Poisson--Nernst--Planck systems~\cite{liu2021positivity, he2019positivity, pan2024linear}.

To our best knowledge, there are very few work on developing thermodynamical consistent two-phase non-Newtonian visco-elastic flow model.  And there is no reported work of structure preserving numerical method for two-phase non-Newtonian visco-elastic flow model. This is the motivation of this work.

In this paper, we first derive the thermodynamical consistent two-phase PENE-P visco-elastic flow model based on phase field method and the variational approach. Then we show that the derived diffusive interface model asymptotically converges to the sharp interface model when the thickness of the diffusive interface model goes to zero.  And this is achieved by performing  singular perturbation asymptotic analysis using interface thickness of the diffusive interface model as the small parameter. And we develop an unconditional energy stable Fourier-Spectral method for solving the  diffusive interface model, where first-order time stepping  is adopted.   Finally, numerical experiments are carried out to confirm the  unconditionally energy stability and convergence order. And the numerical algorithms have been used to simulate droplet motions in the two-phase viscoelastic Couette flow, where the  combination of Fourier and Legendre spectral collocation methods is used.

The rest of this paper is structured as follows. In section~\ref{sec2}, the two-phase FENE-P model based on phase field method is derived by using the thermodynamical consistent framwork. In section~\ref{sec3}, the sharp interface limit of the diffusive interface model is performed. The equivalent system is obtained in section~\ref{sec4}. The unconditionally energy stable scheme is provided in section~\ref{sec5}. Numerical experiments are given in section~\ref{sec6}. And conclusions are given in section~\ref{sec7}.

\section{Two-phase FENE-P viscoelastic fluid model}\label{sec2}
In this section, by adopting the energy variational approach (EnVarA),  we will develop a two-phase FENE-P  viscoelastic fluid model based on the phase field method. 
To be specific, we use a label function to separate two phases with different natures. The interface is described as a smooth transition such that all the variables are smooth functions. 
Besides, the system satisfies the energy decay law.  

\subsection{Sharp interface model}\label{sec:sharp_model}
We now give the sharp interface model of the two-phase FENE-P system 
before we derive the phase field description, to ensure our phase field model 
is based on the reference system. 
We assume the computational domain is $\Omega\subset\mathbb{R}^{\mathrm{dim}},~\mathrm{dim}=2,3$. 
Because of the two-phase FENE-P system, we divide this computational domain 
into two parts represented by $\Omega=\Omega^{+} \bigcup \Omega^{-}$, 
where $\Omega^{+}$ and $\Omega^{-}$ indicate the inner and outer regions respectively, separated by the sharp interface $\Gamma$. Throughout the Section~\ref{sec2}-\ref{sec4}, the convention $(\nabla \bm{u})_{ij}=\partial_ju_i$, $(\nabla \cdot \bm{T} )_i=\sum_{j}\partial_jT_{ij} $ are used, where $\bm{T}$ is the tensor and $\bm{u}$ is the vector. 
The two-phase FENE-P sharp interface model is presented in the following form 
\begin{subequations}
\begin{numcases}{} 
\rho \frac{\mathrm{D}\bm{u}}{\mathrm{D}t} 
+ \nabla p 
= \nabla \cdot \bm{\tau}_{s} 
+ \nabla \cdot \bm{\tau}_{p}, 
& in $\Omega^{\pm}$, 
\label{eqn:sharpU}
\\
\nabla \cdot \bm{u} = 0, 
& in $\Omega^{\pm}$, 
\label{eqn:sharpNablaU}
\\ 
\bm{\tau}_{s} = \mu_{s} \left( \nabla \bm{u} + \nabla \bm{u}^{T}\right)=2  \mu_{s} D \bm{u}, 
& in $\Omega^{\pm}$, 
\label{eqn:sharpTauEta}
\\ 
\bm{\tau}_{p} = \frac{\mu_{p}}{\lambda} f\left( \bm{\sigma} \right) , 
& in $\Omega^{\pm}$, 
\label{eqn:sharpTauP}
\\ 
\frac{\mathrm{D}\bm{\sigma}}{\mathrm{D}t} 
- \nabla \bm{u} \cdot \bm{\sigma}
- \bm{\sigma} \cdot \nabla \bm{u}^{T}
= -\frac{1}{\lambda} f\left( \bm{\sigma} \right), 
& in $\Omega^{\pm}$, 
\label{eqn:sharpSigma}
\\ 
f\left( \bm{\sigma} \right) 
= \frac{\bm{\sigma}}{1-\frac{\mathrm{tr}\left( \bm{\sigma} \right)}{L^{2}}} - \bm{I}, 
& in $\Omega^{\pm}$,
\label{eqn:sharpF}
\end{numcases}
\end{subequations}
where $D\bm{u}=\left(\nabla \bm{u} 
+ \nabla \bm{u}^{T}\right)/2$. Eq. \eqref{eqn:sharpU} and \eqref{eqn:sharpNablaU} is the momentum and mass conservation law for the non-Newtonian viscoelastic FENE-P fluids in two bulk regions, 
where we introduce a common velocity function 
$\bm{u}:\Omega \to \mathbb{R}^{\mathrm{dim}}$ to denote the movements of the non-Newtonian fluid and a general density of the fluid is defined as $\rho$. $\frac{\mathrm{D}}{\mathrm{D}t} = \frac{\partial}{\partial t} + \bm{u} \cdot \nabla$ is  the material derivative. And the densities are two constants in $\Omega^{\pm}$, which means Eq. \eqref{eqn:sharpNablaU} are well established. 
Expression $\bm{\tau}_{s}$ and $\bm{\tau}_{p}$ are the viscous stress and 
elastic (polymeric) stress respectively, where $\mu_{s}$ is the solvent viscosity, 
$\mu_{p}$ and $\lambda$ define the polymer viscosity and relaxation time, 
respectively. $p$ is the pressure, and $\bm{\sigma}$ is the conformation tensor which is dimensionless quantity and positive definite, 
$\mbox{tr}\left(\cdot\right)$ denotes the trace of a matrix tensor, 
$L$ represents the dimensionless finite extensibility length. 
Eq. \eqref{eqn:sharpSigma} gives the movement description 
for the conformation tensor $\bm{\sigma}$ of the polymer. 
$f$ is the relaxation function, $\bm{I}$ is the unitary tensor. 


We should specify two jump conditions at the two-phase fluid interface $\Gamma$, 
\begin{equation}\label{eqn:sharp_condition}
\left[\bm{u}\right]_{-}^{+} = 0, \quad 
\left[ - p \bm{I}+ \bm{\tau}_{s} 
+\bm{\tau}_{p} 
 \right]_{-}^{+} \cdot \bm{n}
= -T \kappa \bm{n}, 
\end{equation} 
where $\left[ f \right]_{-}^{+}=f^{+}-f^{-}$ denotes the jump of function $f$ across the interface, $\kappa=-\nabla_s\cdot \bm{n}$ is the mean curvature of the interface,  $T$ is the interfacial tension coefficient between two fluids, 
and $\bm{n}$ is the unit normal vector on $\Gamma$ pointing from $\Omega_{-}$ to $\Omega_{+}$. 
The first condition of Eq. \eqref{eqn:sharp_condition} is the continuity of the velocity across the interface. 
The second condition of Eq. \eqref{eqn:sharp_condition} is the force balance across the interface including the total stress and the interfacial tension  at the sharp interface.

\subsection{Phase field modeling}
In this section, we derive the two-phase FENE-P phase field model based on energy variational approach method \cite{Qin2022JCP}. 
We first introduce an order parameter $\varphi$ to describe two separated phases with a transition in the smooth interface 
distinguished by $\varphi=\pm 1$, with a computational domain $\Omega\subset\mathbb{R}^{\mathrm{dim}},~\mathrm{dim}=2,3$, which has the similar definition shown in section \ref{sec:sharp_model}. 
The inner and outer regions $\Omega^{\pm}$ are approximated by the domain occupied by two different phases, 
separated by a diffuse interface formulated by $\varphi \in \left(-1,1\right)$, 
which is different from section \ref{sec:sharp_model}. 
We define the velocity field to be the volume-averaged velocity of this mixture, 
which is given a detailed description in \cite{Abels2012Thermodynamically}. 
We choose the density $\rho$ as a linear combination of the densities in the two phases: 
$\rho\left(\varphi\right)=\frac{1+\varphi}{2}\rho^{+}+\frac{1-\varphi}{2}\rho^{-}$. 
Based on the assumption of mass conservation and momentum conservation, 
We give the following kinetic equations for label function $\varphi$ and velocity $\bm{u}$. 
\begin{subequations}\label{eqn:main}
\begin{numcases}
{}
\frac{\mathrm{D}\varphi}{\mathrm{D}t} + \nabla \cdot \bm{J} = 0, 
& in $\Omega$, 
\label{eqn:DimPhi}\\
\frac{\mathrm{D}\left(\rho\bm{u}\right)}{\mathrm{D}t} + \nabla p 
= \nabla \cdot \bm{\tau}_{s} 
+ \nabla \cdot \bm{\tau}_{p} 
+ \nabla \cdot \bm{\tau}_{\varphi}, 
& in $\Omega$, 
\label{eqn:DimU}\\
\nabla \cdot \bm{u} = 0, 
& in $\Omega$, 
\label{eqn:DimNablaU}
\end{numcases}
\end{subequations}
where the viscosity stress 
$\bm{\tau}_{s} = \mu_{s} \left( \varphi \right) \left( \nabla \bm{u} + \nabla \bm{u}^{T}\right)$,  the solvent viscosity $ \mu_{s} $ is chosen as a linear combination of solvent viscosities in the two phases: $\mu_s\left(\varphi\right)=\frac{1+\varphi}{2}\mu^{+}_s+\frac{1-\varphi}{2}\mu^{-}_s$, and the elastic stress $\bm{\tau}_{p}$ also depends on 
the phase field function $\varphi$. 
$\bm{\tau}_{\varphi}$ is the surface tension force. 
The flux $\bm{J}$ is specified later to meet the requirement of non-increasing energy. 
$p$ is the pressure and operator $\frac{\mathrm{D}}{\mathrm{D}t} = \frac{\partial}{\partial t} + \bm{u} \cdot \nabla$ 
is the material derivative. Equation \eqref{eqn:DimNablaU} 
is the impressible condition derived from the mass conservation. 
To obtain an equation for the diffuse elastic strain tensor, 
we need to combine the variables in two fluid phases by phase field function.  To do this, we first give the equation for  the conformation tensor $\bm{\sigma}$
\begin{equation}
\frac{\mathrm{D}\bm{\sigma}}{\mathrm{D}t} 
- \nabla \bm{u} \cdot \bm{\sigma}
- \bm{\sigma} \cdot \nabla \bm{u}^{T}
= -\frac{1}{\lambda\left(\varphi\right)}f \left( \bm{\sigma} \right),\quad f\left( \bm{\sigma} \right) 
= \frac{\bm{\sigma}}{1-\frac{\mathrm{tr}\left( \bm{\sigma} \right)}{L^{2}}} - \bm{I},
\end{equation}
where 
$\lambda \left( \varphi \right)$ is the relaxation time also chosen as linear combination of the constant quantities from two fluid phases: $\lambda \left( \varphi \right)=\frac{1+\varphi}{2}\lambda^{+}+\frac{1-\varphi}{2}\lambda^{-}$. 

The total phase-dependent elastic stress is then given by the following formulation: 
\begin{equation} 
\bm{\tau}_{p} 
= \frac{\mu_{p}\left( \varphi \right)}{\lambda\left( \varphi \right)} 
f \left( \bm{\sigma} \right), 
\end{equation}
where $\mu_{p}\left( \varphi \right)$ is the polymer viscosity also chosen as linear combination of the constant quantities from two fluid phases: $\mu_{p} \left( \varphi \right)=\frac{1+\varphi}{2}\mu^{+}_{p}+\frac{1-\varphi}{2}\mu^{-}_{p}$. 

Following the energy definition in \cite{Hu2007New}, the total energy of the two-phase flow system is defined as 
\begin{equation}
\begin{aligned}
& E_{total} = E_{kin} + E_{ela} + E_{mix} 
\\ 
& = \underbrace{\int_{\Omega} \frac{\rho\left(\varphi\right)}{2} \left| \bm{u} \right|^{2}\mathrm{d}\bm{x}}_{\mbox{kinetic energy}} 
\underbrace{- \int_{\Omega}\frac{\mu_{p}\left(\varphi\right)}{2\lambda\left(\varphi\right)} 
\left( L^{2} \ln \left( 1-\frac{\mathrm{tr}\left(\bm{\sigma}\right)}{L^{2}}\right) 
+ \mathrm{tr} \left(\ln\left(\bm{\sigma}\right) + \bm{I} \right) \right)\mathrm{d}\bm{x}}_{\mbox{elastic energy}} 
\underbrace{+ \int_{\Omega}\gamma \left( \frac{\delta}{2} \left| \nabla \varphi \right|^{2} 
+ \frac{1}{\delta} W \left( \varphi \right) \right) \mathrm{d} \bm{x}}_{\mbox{mixing energy}}, 
\end{aligned}
\end{equation}
where $\gamma$ is the surface tension density, 
$\delta$ is the thickness of the diffuse interface 
and $W\left(\varphi\right)$ is the double-well potential to describe the two phases, which is chosen as $W\left(\varphi\right) = \left( 1- \varphi^{2} \right)^{2}/4$ throughout this work. 

For simplicity, we consider a closed system in $\Omega$ without the energy input, which means the boundary effects are neglected from integration parts based on some appropriate boundary conditions.
 Specifically, the boundary conditions are periodic boundary condition or the following zero Neumann condition for phase field variable $\varphi$, zero flux condition and zero velocity condition:
\begin{equation} 
\left. \nabla \varphi \cdot \bm{n} \right|_{\partial \Omega} = 0, 
\quad 
\left. \bm{J} \cdot \bm{n} \right|_{\partial \Omega} = 0, 
\quad 
\left. \bm{u} \right|_{\partial \Omega} = 0,  
\end{equation}\label{BC_zero} 
where $\bm{n}$ is the unit outer normal vector on the boundary. 

 By the energy variational approach, we will derive the following  FENE-P-CH  phase field model:
\begin{subequations}\label{eqn:system}
\begin{numcases}{} 
\frac{\mathrm{D}\varphi}{\mathrm{D}t} = -\nabla \cdot \bm{J} 
= \nabla \cdot \left( M \left( \varphi \right) \nabla \nu \right), \quad \bm{J}=-M \left( \varphi \right) \nabla \nu,
& in $\Omega$, 
\label{eqn:mainPhi}
\\
\nu = \gamma \left( - \delta \Delta \varphi 
+ \frac{1}{\delta} W^{\prime} \left( \varphi \right) \right) 
- \left( \frac{\mu_{p}\left( \varphi \right)}{2\lambda\left( \varphi \right)} \right)^{\prime}
\left( L^{2} \ln \left(1-\frac{\mathrm{tr}\left(\bm{\sigma}\right)}{L^{2}}\right) 
+ \mathrm{tr} \left(\ln \left( \bm{\sigma}\right)+\bm{I} \right) \right), 
& in $\Omega$, 
\label{eqn:mainMu}
\\ 
\frac{\mathrm{D}\left(\rho\bm{u}\right)}{\mathrm{D}t} 
-\nabla \cdot \left( \rho^{\prime} \left( \varphi \right) \bm{u} \otimes M\left( \varphi \right) \nabla \nu \right) + \nabla p 
= \nabla \cdot \bm{\tau}_{s} 
+ \nabla \cdot \bm{\tau}_{p} 
+ \nabla \cdot \bm{\hat{\tau}}_{\varphi}, 
& in $\Omega$, 
\label{eqn:mainU}
\\
\nabla \cdot \bm{u} = 0, 
& in $\Omega$, 
\label{eqn:mainNablaU}
\\ 
\bm{\tau}_{s} = \mu_{s} \left( \varphi \right) \left( \nabla \bm{u} + \nabla \bm{u}^{T}\right), 
& in $\Omega$, 
\label{eqn:mainTauEta}
\\ 
\bm{\tau}_{p} = \frac{\mu_{p}\left( \varphi \right)}{\lambda\left( \varphi \right)} f(\bm {\sigma}), 
& in $\Omega$, 
\label{eqn:mainTauP}
\\ 
\bm{\hat{\tau}}_{\varphi} = - \gamma \delta \nabla \varphi \otimes \nabla \varphi, 
& in $\Omega$, 
\label{eqn:mainTauPhi} 
\\ 
\frac{\mathrm{D}\bm{\sigma}}{\mathrm{D}t} 
- \nabla \bm{u} \cdot \bm{\sigma}
- \bm{\sigma} \cdot \nabla \bm{u}^{T}
= -\frac{1}{\lambda\left(\varphi\right)}f \left( \bm{\sigma} \right), 
& in $\Omega$, 
\label{eqn:mainSigma}
\\ 
f\left( \bm{\sigma} \right) 
= \frac{\bm{\sigma}}{1-\frac{\mathrm{tr}\left( \bm{\sigma} \right)}{L^{2}}} - \bm{I}, 
& in $\Omega$. 
\label{eqn:mainF}
\end{numcases}
\end{subequations}
where the symbol $^{\prime}$ denotes the derivative of $\varphi$. 
And the above model satisfies the following energy dissipative relation:
\begin{equation}\label{eqn:dissipation}
\begin{aligned}
\frac{\mathrm{d}E_{total}}{\mathrm{d}t} 
= & \frac{\mathrm{d}E_{kin}}{\mathrm{d}t} 
+ \frac{\mathrm{d}E_{ela}}{\mathrm{d}t} 
+ \frac{\mathrm{d}E_{mix}}{\mathrm{d}t} 
\\ 
= & - \int_{\Omega} \frac{\mu_{s} \left( \varphi \right)}{2} 
\left| \nabla \bm{u} + \nabla \bm{u}^{T}\right|^{2} \mathrm{d}\bm{x} - \int_{\Omega} M(\phi)|\nabla \nu|^2\bm{x} - \int_{\Omega} \frac{\mu_{p}\left( \varphi \right)}{2\lambda^{2}\left( \varphi \right)} 
\mathrm{tr}\left(f^{2} \left( \bm{\sigma} \right) \bm{\sigma}^{-1} \right) 
\mathrm{d} \bm{x}. 
\end{aligned}
\end{equation}
The detailed derivation of the above FENE-P-CH model is listed in the appendix~\ref{labelapp2}.


\subsection{Dimensionless governing equation}\label{sec:dimensionlization}
To make the analysis more simplifier, 
we provide the dimensionless governing equations in this section. 
We choose the characteristic length $\mathbb{L}$ and characteristic velocity $\mathbb{U}$ so that the characteristic time can be written as $\mathbb{T} = \mathbb{L} / \mathbb{U}$. 
As we consider the two-phase FENE-P model, 
there are four pairs of coefficients, 
they are solvent viscosities $\mu^{+}_s$ and $\mu^{-}_s$, densities $\rho_{+}$ and $\rho_{-}$, 
polymer viscosities $\mu^{+}_{p}$ and $\mu^{-}_{p}$ 
and relaxation times $\lambda^{+}$ and $\lambda^{-}$. 
Each pair of them is connected by the label function $\varphi$ using algebraic average. 
Here we choose $\rho$, $\mu_s$, $\mu_{p}$ and $\lambda$ as the characteristic density, characteristic solvent viscosity, characteristic polymer viscosity and relaxation time.
We also choose the characteristic mobility as $\mathbb{M}$. 
The characteristic chemical potential is defined as $\gamma/\mathbb{L}$. The characteristic values of  pressure, the viscous stress $\tau_{s}$, the polymer stress $\tau_p$ and surface stress $\bm{\hat{\tau}}_{\varphi}$ are all chosen as $\mu_{s}\mathbb{U}/\mathbb{L}$. 
All the dimensional variables can be written as 
\begin{equation*} 
\begin{aligned} 
& x = \tilde{x} \mathbb{L}, \quad 
\bm{u} = \tilde{\bm{u}} \mathbb{U}, \quad 
t = \tilde{t} \mathbb{T}, \quad 
M \left( \varphi \right) = \tilde{M} \left( \varphi \right) \mathbb{M}, \quad 
\delta = \tilde{\delta} \mathbb{L}, \quad 
\nu = \tilde{\nu} \frac{\gamma}{\mathbb{L}},\quad 
p=\tilde{p}\frac{\mu_{s}\mathbb{U}}{\mathbb{L}},\\ 
& \mu_{s} \left( \varphi \right) 
= \left( \frac{1-\varphi}{2}\tilde{\mu}^{-}_{s} + \frac{1+\varphi}{2}\tilde{\mu}^{+}_{s} \right) \mu_{s} 
= \tilde{\mu}_{s} \left( \varphi \right) \mu_{s}, \quad 
\rho \left( \varphi \right) 
= \left( \frac{1-\varphi}{2}\tilde{\rho}^{-} + \frac{1+\varphi}{2}\tilde{\rho}^{+}\right) \rho
= \tilde{\rho} \left( \varphi \right) \rho, \\ 
&\mu_{p} \left( \varphi \right) 
= \left( \frac{1-\varphi}{2} \tilde{\mu}^{-}_{p}  + \frac{1+\varphi}{2} \tilde{\mu}^{-}_{p} \right) \mu_{p} 
= \tilde{\mu}_{p} \left( \varphi \right) \mu_{p}, \quad 
\lambda \left( \varphi \right) 
= \left( \frac{1-\varphi}{2}\tilde{\lambda}^{-} + \frac{1+\varphi}{2} \tilde{\lambda}^{+} \right) \lambda
= \tilde{\lambda} \left( \varphi \right) \lambda,
\end{aligned}
\end{equation*} 
where
\[\tilde{\mu}^{\pm}_{s}=\frac{\mu^{\pm}_{s}}{\mu_s},\quad \tilde{\rho}^{\pm}=\frac{\rho^{\pm}}{\rho},\quad 
\tilde{\mu}^{\pm}_{p}=\frac{\mu^{\pm}_{p}}{\mu_p},\quad \tilde{\lambda}^{\pm}=\frac{\lambda^{\pm}}{\lambda}.\]
All operators have  similar definitions in which the dimensionless operators have a tilde symbol over the head, and we omit the complicated formulations. 
After careful calculations and omitting the symbol $\tilde{\cdot}$  for all dimensionless variables and operators without ambiguity, the dimensionless FENE-P-CH system can be written as 

\begin{subequations}\label{eqn:main_system}
\begin{numcases}{} 
\frac{\mathrm{D}\varphi}{\mathrm{D}t} 
= \frac{1}{\mathrm{Pe}} \nabla \cdot \left( M \left( \varphi \right) 
\nabla \nu \right), 
& in $\Omega$, 
\label{eqn:phi}
\\
\nu = - \delta \Delta \varphi 
+ \frac{1}{\delta} W^{\prime} \left( \varphi \right) 
- \frac{\beta\mathrm{Ca}}{\mathrm{Wi}} \left( \frac{\mu_{p}\left( \varphi \right)}{2\lambda\left( \varphi \right)} \right)^{\prime}
\left( L^{2} \ln \left(1-\frac{\mathrm{tr}\left(\bm{\sigma}\right)}{L^{2}}\right) 
+ \mathrm{tr} \left(\ln \left( \bm{\sigma}\right)+\bm{I} \right) \right), 
& in $\Omega$, 
\label{eqn:mu}
\\ 
\mathrm{Re} \left( \frac{\mathrm{D}\left(\rho\bm{u}\right)}{\mathrm{D}t} 
- \frac{1}{\mathrm{Pe}} \nabla \cdot \left( \rho^{\prime} \left( \varphi \right) \bm{u} \otimes M \left( \varphi \right) \nabla \nu \right) \right) + \nabla p 
= \nabla \cdot \bm{\tau}_{s} 
+ \nabla \cdot \bm{\tau}_{p} 
+\nabla \cdot \bm{\tau}_{\varphi}, 
& in $\Omega$, 
\label{eqn:u}
\\
\nabla \cdot \bm{u} = 0, 
& in $\Omega$, 
\label{eqn:nablau}
\\ 
\bm{\tau}_{s} = \mu_{s} \left( \varphi \right) \left( \nabla \bm{u} + \nabla \bm{u}^{T}\right), 
& in $\Omega$, 
\label{eqn:taueta}
\\ 
\bm{\tau}_{p} = \frac{\beta}{\mathrm{Wi}}\frac{\mu_{p}\left( \varphi \right)}{\lambda\left( \varphi \right)} 
f(\bm {\sigma}), 
& in $\Omega$, 
\label{eqn:taup}
\\ 
\bm{\tau}_{\varphi} = - \frac{\delta}{\mathrm{Ca}}  \nabla \varphi \otimes \nabla \varphi, 
& in $\Omega$, 
\label{eqn:tauphi} 
\\ 
\frac{\mathrm{D}\bm{\sigma}}{\mathrm{D}t} 
- \nabla \bm{u} \cdot \bm{\sigma}
- \bm{\sigma} \cdot \nabla \bm{u}^{T}
= -\frac{1}{\mathrm{Wi}\lambda\left(\varphi\right)}f \left( \bm{\sigma} \right), 
& in $\Omega$, 
\label{eqn:sigma}
\\ 
f\left( \bm{\sigma} \right) 
= \frac{\bm{\sigma}}{1-\frac{\mathrm{tr}\left( \bm{\sigma} \right)}{L^{2}}} - \bm{I}, 
& in $\Omega$, 
\label{eqn:f}
\end{numcases}
\end{subequations}
where
\begin{equation*} 
\begin{aligned} 
\mu_{s} \left( \varphi \right)&= \frac{1-\varphi}{2}\mu^{-}_s + \frac{1+\varphi}{2}\mu^{+}_s, \quad
\mu_{p} \left( \varphi \right)=\frac{1-\varphi}{2}\mu^{-}_{p} + \frac{1+\varphi}{2}\mu^{+}_{p},  \\
\rho \left( \varphi \right)&=\frac{1-\varphi}{2}\rho^{-} + \frac{1+\varphi}{2}\rho^{+} ,\quad
\lambda \left( \varphi \right)= \frac{1-\varphi}{2}\lambda^{-} + \frac{1+\varphi}{2}\lambda^{+}, \\
\end{aligned} 
\end{equation*} 
and $\beta$ is a dimensionless number to measure the ratio of polymer viscosity and solvent viscosity in $\Omega_{-}$, 
and the dimensionless Peclet number $\mathrm{Pe}$, Reynolds number $\mathrm{Re}$, 
Capillary number $\mathrm{Ca}$ and Weissenberg number $\mathrm{Wi}$ are given as follows 
\begin{equation}
\begin{aligned}
\mathrm{Pe} = \frac{\mathbb{L}^{3}}{\gamma \mathbb{M} \mathbb{T}}, \quad 
\mathrm{Re} = \frac{\mathbb{U} \mathbb{L}\rho}{\mu_{s}}, \quad 
\mathrm{Ca} = \frac{\mathbb{U} \mu_{s}}{\gamma}, \quad 
\mathrm{Wi} = \frac{\mathbb{U} \lambda}{\mathbb{L}}, \quad 
\beta = \frac{\mu_{p}}{\mu_{s}}. 
\end{aligned}
\end{equation}

The initial conditions are given by 
\[\bm{u}(\bm{x},0)=\bm{u}_0,\quad \phi(\bm{x},0)=\phi_0,\quad \bm{\sigma}(\bm{x},0)=\bm{\sigma}_0. \]

System \eqref{eqn:main_system} satisfies the following energy dissipative law. 
\begin{theorem}[Energy dissipation]\label{thm:energy}
    System \eqref{eqn:main_system} satisfies the following dissipative property: 
    \begin{equation}
\begin{aligned}
\frac{\mathrm{d}E_{total}}{\mathrm{d}t} 
= & \frac{\mathrm{d}E_{kin}}{\mathrm{d}t} 
+ \frac{\mathrm{d}E_{ela}}{\mathrm{d}t} 
+ \frac{\mathrm{d}E_{mix}}{\mathrm{d}t} 
\\ 
= & - \int_{\Omega} \frac{\mu_{s} \left( \varphi \right)}{2} 
\left| \nabla \bm{u} + \nabla \bm{u}^{T}\right|^{2} \mathrm{d}\bm{x} 
- \int_{\Omega} \frac{M(\phi)}{\mathrm{CaPe}}\left|\nabla \nu\right|^2\bm{x} 
- \int_{\Omega} \frac{\beta\mu_{p}\left( \varphi \right)}{2\mathrm{Wi}^{2}\lambda^{2}\left( \varphi \right)} 
\mathrm{tr}\left(f^{2} \left( \bm{\sigma} \right) \bm{\sigma}^{-1} \right) 
\mathrm{d} \bm{x}, 
\end{aligned}
\end{equation}
where the total energy is formulated as 
\begin{equation}
\begin{aligned}
    & E_{total} 
    = E_{kin} + E_{mix} + E_{ela} 
    \\ 
    & = \int_{\Omega} 
    \underbrace{\frac{\mathrm{Re}}{2} \rho \left| \bm{u} \right|^{2}}_{kinetic \ energy} 
    +  \underbrace{\frac{1}{\mathrm{Ca}} \left( \frac{\delta}{2}\left|\nabla\varphi\right|^{2} 
    + \frac{1}{\delta} W \left( \varphi \right) \right)}_{mixed \ energy} 
     \underbrace{-\frac{\beta}{2\mathrm{Wi}}\frac{\mu_{p}\left( \varphi \right)}{\lambda\left( \varphi \right)}
\left( L^{2} \ln \left( 1-\frac{\mathrm{tr}\left(\bm{\sigma}\right)}{L^{2}}\right) 
+ \mathrm{tr} \left(\ln\left(\bm{\sigma}\right) + \bm{I} \right) \right)}_{elastic \ energy}
\mathrm{d}\bm{x}. 
\end{aligned}
\end{equation}
\end{theorem}
In the rest of the work, we assume mobility $M\left( \varphi \right)$ is a constant, which is order of $\delta $, i.e.,
     $M\left( \varphi \right) = M_0\delta $, where $M_0$ is a constant. 

\section{Sharp interface limit analysis}\label{sec3}
We model the FENE-P-CH system based on the phase field approach with a diffusive interface introduced. 
Therefore, the sharp interface limit when the thickness $\delta\to 0$ needs to be considered. 
In this section, we talk about the sharp interface limit for the FENE-P-CH model. 
We follow the basic idea in \cite{Abels2012Thermodynamically} to perform the formal matched asymptotic expansions. 
Therefore, we need to consider the solution $\left( \bm{u}, p, \rho, \varphi, \nu, \bm{\sigma} \right)$ of the system given in \eqref{eqn:main_system}. 
\subsection{Outer expansions} 
We need to expand the solution in the bulk regions far away the diffuse interface. 
We assume the expansion is of the following form 
\begin{equation}\label{exp:out}
\begin{aligned}
    & \bm{u}^{\pm} = \sum_{k=0}^{\infty} \delta^{k} \bm{u}_{k}^{\pm}, \quad 
     p^{\pm} = \sum_{k=-1}^{\infty} \delta^{k} p_{k}^{\pm}, \quad 
    \bm{\tau}_{\eta}^{\pm} 
    = \sum_{k=0}^{\infty} \delta^{k} \bm{\tau}_{\eta,k}^{\pm}, \quad 
    \bm{\tau}_{p}^{\pm} 
    = \sum_{k=0}^{\infty} \delta^{k} \bm{\tau}_{p,k}^{\pm}, \\ 
    & \varphi^{\pm} = \sum_{k=0}^{\infty} \delta^{k} \varphi_{k}^{\pm}, \quad 
    \nu^{\pm} = \sum_{k=0}^{\infty} \delta^{k} \nu_{k}^{\pm}, \quad 
    \bm{\sigma}^{\pm} = \sum_{k=0}^{\infty} \delta^{k} \bm{\sigma}_{k}^{\pm},\quad\bm{\tau}_{\varphi}^{\pm} 
    = \sum_{k=0}^{\infty} \delta^{k} \bm{\tau}_{\varphi,k}^{\pm}. 
\end{aligned}
\end{equation}
The expansion of Eq. \eqref{eqn:mu} at the order $\delta^{-1}$ leads to 
$W^{\prime}\left(\varphi_{0}^{\pm}\right)=0$, 
which yields the stable solutions $\varphi_{0}^{\pm} = \pm 1$. We denote by $\Omega^{\pm}$ the regions where $\varphi^{\pm}_{0}= \pm1$. The expansion of Eq.\eqref{eqn:u}  at the order $\delta^{-1}$ leads to 
$\nabla p^{\pm}_{-1}=0$, which gives $p^{\pm}_{-1}=c^{\pm}$, where $c^{\pm}$ are two constants in space.

Insert the expansions \eqref{exp:out} into system \eqref{eqn:main_system}, as $M\left( \varphi \right) = M_0\delta=O(\delta)$, at the order of O($1$),  we obtain 
\begin{subequations}
\begin{numcases}{} 
\mathrm{Re} \rho(\varphi^{\pm}_0) \frac{\mathrm{D}\bm{u}^{\pm}_0}{\mathrm{D}t} 
+ \nabla p^{\pm}_0 
= \nabla \cdot \bm{\tau}_{s,0}^{\pm} 
+ \nabla \cdot \bm{\tau}_{p,0}^{\pm}, 
& in $\Omega_{\pm}$, 
\label{eqn:outerU}
\\
\nabla \cdot \bm{u}^{\pm}_0 = 0, 
& in $\Omega_{\pm}$, 
\label{eqn:outerNablaU}
\\ 
\bm{\tau}_{\eta,0}^{\pm} = \mu_{s}(\varphi^{\pm}_0) \left( \nabla \bm{u}^{\pm}_0 + \left(\nabla \bm{u}^{\pm}_0\right)^{T}\right), 
& in $\Omega_{\pm}$, 
\label{eqn:outerTauEta}
\\ 
\bm{\tau}_{p,0}^{\pm} = \frac{\beta}{\mathrm{Wi}} \frac{\mu_{p}(\varphi^{\pm})}{\lambda(\varphi^{\pm}_0)} 
f(\bm {\sigma^{\pm}_0}), 
& in $\Omega_{\pm}$, 
\label{eqn:outerTauP}
\\ 
\bm{\tau}_{\varphi,0}^{\pm}=0, & in $\Omega_{\pm}$, 
\label{eqn:outerTauPhi}
\\ 
\frac{\mathrm{D}\bm{\sigma}^{\pm}_0}{\mathrm{D}t} 
- \nabla \bm{u}^{\pm}_0 \cdot \bm{\sigma}^{\pm}_0
- \bm{\sigma}^{\pm}_0 \cdot \left(\nabla \bm{u}^{\pm}_0\right)^{T}
= -\frac{1}{\mathrm{Wi}\lambda(\varphi^{\pm}_0)} f\left( \bm{\sigma}^{\pm}_0 \right), 
& in $\Omega_{\pm}$, 
\label{eqn:outerSigma}
\\ 
f\left( \bm{\sigma}^{\pm}_0 \right) 
= \frac{\bm{\sigma}^{\pm}_0}{1-\frac{\mathrm{tr}\left( \bm{\sigma}_0^{\pm} \right)}{L^{2}}} - \bm{I}, 
& in $\Omega_{\pm}$,  
\label{eqn:outerF}
\end{numcases}
\end{subequations}
where 
\begin{equation*}
\rho(\varphi^{\pm}_0)=\rho^{\pm},\quad \mu_s(\varphi^{\pm}_0)=\mu^{\pm}_{s},\quad \mu_p(\varphi^{\pm}_0)=\mu^{\pm}_{p},\quad \lambda(\varphi^{\pm}_0)=\lambda^{\pm}.
\end{equation*}
The above leading equations recover the sharp interface equations given in \eqref{eqn:sharpU}-\eqref{eqn:sharpF}. 
\begin{remark}
    Noting that there are three different roots for the equation $W^{\prime}\left(\varphi_0\right)=0$, namely $\varphi_0=-1,0,1$. In this subsection we consider the outer region, which means the bulk regions are away from the interface. Therefore $\varphi_{0}=0$ is not considered in the outer expansions. 
\end{remark}
\subsection{Coordinate transformation and inner expansions} 
In this section, we use the similar concepts shown in \cite{Abels2012Thermodynamically} to make the analysis more easier. 
For simplicity, we restate the main results and definitions again. 
We denote the smoothly evolving interface as 
$\Gamma = \Gamma\left( t \right)_{t \geq 0}$, 
which we expect to be the limit of the zero-level sets of the order parameter 
$\varphi$ when the thickness of diffuse interface $\delta$ tends to zero, 
and the zero limit interface is represented as $\Gamma^{0}\left( t \right)$. 
Firstly, the expansion in the diffuse interface region is given 
and therefore new coordinates are well defined in a neighborhood of the interface. 
We introduce the following local parameterization of $\Gamma$  
\begin{equation}
    \bm{\zeta}: I \times U \to \mathbb{R}^{d}, 
\end{equation}
where the time interval $I\in\mathbb{R}$ 
and spatial domain $U\subset\mathbb{R}^{d-1}$ are chosen. 
The unit normal to $\Gamma\left(t\right)$ will be denoted by $\bm{n}$, 
which is the same as in the previous section \ref{sec:sharp_model} 
and points from $\Omega_{-}$ into $\Omega_{+}$. 
The signed distance function $g\left(t,\bm{x}\right)$ of a point $\bm{x}$ is introduced to the interface $\Gamma^{0}\left( t \right)$ with $g\left( t , \bm{x} \right) > 0 $ if $\bm{x} \in \Omega_{+}$. 
Based on the above definition, 
and by using a rescaled distance $z = g \left( t , \bm{x} \right) / \delta$,  we now introduce the new coordinates 
by defining a local parameterization of $I\times \mathbb{R}^{d}$ 
close to $\bm{\zeta}\left(I\times U\right)$ as follows: 
\begin{equation}
    G^{\delta} \left( t, \bm{s}, z \right) 
    = \left( t, \bm{\zeta} \left( t, \bm{s} \right) 
    + \delta z \bm{n} \left( t, \bm{s} \right) \right) 
\end{equation}
with $\bm{s} \in U \subset \mathbb{R}^{d-1}$. 
We define 
\begin{equation}
    V= \partial_{t} \bm{\zeta} \cdot \bm{n} 
\end{equation}
as the (scalar) normal velocity 
and  the inverse function  $(t,\bm{s},z)(t,\bm{x}):=(G^{\delta})^{-1}(t,\bm{x})$ satisfies:
\begin{equation}
    \partial_{t} z = \frac{1}{\delta} \partial_{t} g 
    = -\frac{1}{\delta} V, 
\end{equation}
which is a result of \cite{Abels2012Thermodynamically} 
and the theoretical proof can be concluded based on (2.6) and (2.20) in 
\cite{Deckelnick2005Computation}. 
Based on the coordinate transformation, we can rewrite the variables 
with the new coordinates in the following forms containing a scalar function 
$b\left( t, \bm{x} \right)$ and 
a vector field $\bm{j} \left( t, \bm{x} \right)$: 
\begin{subequations}
\begin{align}
    \hat{b} \left( t,s\left(t,\bm{x}\right), z\left(t,\bm{x}\right)\right)
    &=b\left(t,\bm{x}\right), 
    \\ 
    \hat{\bm{j}} \left( t,s\left(t,\bm{x}\right), z\left(t,\bm{x}\right)\right)
    &=\bm{j}\left(t,\bm{x}\right). 
\end{align}
\end{subequations}
Hence for a scalar function $b\left( t, \bm{x} \right)$ and 
a vector field $\bm{j} \left( t, \bm{x} \right)$, 
we have the following identities to complete the asymptotic analysis, 
\begin{subequations}
    \begin{align}
        \frac{\mathrm{d}}{\mathrm{d}t} b \left( t, \bm{x} \right) 
         &= -\frac{1}{\delta} V \partial_{z} \hat{b}+ \partial_{t}\bm{s}\cdot\nabla_{\bm{s}} \hat{b}+\partial_{t}\hat{b}= -\frac{1}{\delta} V \partial_{z} \hat{b} + \mbox{h.o.t.}, 
         \\ 
         \nabla_{\bm{x}}b 
         &= \nabla_{\Gamma_{\delta z}} \hat{b} + \frac{1}{\delta} \partial_{z} \hat{b}\bm{n}, 
         \label{eqn:inner_grad}
         \\ 
         \nabla_{\bm{x}} \cdot \bm{j}  
         &= \nabla_{\Gamma_{\delta z}} \cdot \hat{\bm{j}} + \frac{1}{\delta} \partial_{z}\hat{\bm{j}}\cdot \bm{n}, 
    \end{align}
\end{subequations}
where we use h.o.t. to stand for higher order terms as the thickness of the interface trends to zero i.e. $\delta \to 0$, and the following correspondence establishes  
\begin{enumerate}
    \item $\nabla_{\bm{x}}$ and $\nabla_{\bm{x}}\cdot$ represent 
    the gradient and divergence concerning the spatial variables $\bm{x}$, 
    respectively.
    \item $\nabla_{\Gamma_{\delta z}}$ is the surface gradient on 
    $\Gamma_{\delta z}:=\left\{\bm{\zeta}\left( \bm{s} \right) 
    + \delta z \bm{n}\left( \bm{s} \right)| s\in U\right\}$. 
\end{enumerate}
For the convenience of discussion, without causing ambiguity, 
the time dependence is ignored in the following argumentations. 
Moreover, we will make use of the relations in our analysis 
\begin{subequations}
    \begin{align}
        \nabla_{\Gamma_{\delta z}} \hat{b} \left( \bm{s}, z\right) 
        & = \nabla_{\Gamma}\hat{b}\left(\bm{s},z\right) + \mbox{h.o.t.}, 
        \\ 
        \nabla_{\Gamma_{\delta z}} \cdot \hat{\bm{j}} \left( \bm{s}, z \right) 
        & = \nabla_{\Gamma} \cdot \hat{\bm{j}} \left( \bm{s},z\right) + \mbox{h.o.t.}, 
        \label{inner_exp:div}
        \\ 
        \Delta_{\Gamma_{\delta z}} \hat{b}\left( \bm{s},z\right) 
        & = \Delta_{\Gamma}\hat{b} \left( \bm{s}, z \right) + \mbox{h.o.t.}, 
    \end{align}
\end{subequations}
where $\nabla_{\Gamma}$, $\nabla_{\Gamma}\cdot$ and $\Delta_{\Gamma}$ 
are the respective surface gradient, divergence, and Laplacian operators on $\Gamma$. The comprehensive theoretical proof of the above identities can be found in the appendix of Ref. \cite{Abels2012Thermodynamically}, 
therefore we omit the details in our paper. 
Similar to the outer variables, we now assume a $\delta$-series approximation for the unknown functions $\left( \bm{u}^{\delta}, p^{\delta}, \rho^{\delta}, \varphi^{\delta}, \nu^{\delta}, \bm{\sigma}^{\delta} \right)$ for the inner variables in the diffuse interface: 
\begin{equation}\label{exp:in}
    \begin{aligned}
        & \bm{u}^{\delta} = \sum_{k=0}^{\infty} \delta^{k} \bm{u}_{k}^{\delta}, \quad 
         p^{\delta} = \sum_{k=-1}^{\infty} \delta^{k} p_{k}^{\delta}, \quad 
        \bm{\tau}_{\eta}^{\delta} = \sum_{k=0}^{\infty} \delta^{k} \bm{\tau}_{\eta,k}^{\delta}, \quad 
        \bm{\tau}_{p}^{\delta} = \sum_{k=0}^{\infty} \delta^{k} \bm{\tau}_{p,k}^{\delta}, \\ 
        & \varphi^{\delta} = \sum_{k=0}^{\infty} \delta^{k} \varphi_{k}^{\delta}, \quad 
        \nu^{\delta} = \sum_{k=0}^{\infty} \delta^{k} \nu_{k}^{\delta}, \quad 
        \bm{\sigma}^{\delta} = \sum_{k=0}^{\infty} \delta^{k} \bm{\sigma}_{k}^{\delta},\quad 
        \bm{\tau}^{\delta}_{\varphi}  = \sum_{k=0}^{\infty} \delta^{k} \bm{\tau}^{\delta}_{\varphi,k}. 
    \end{aligned}
\end{equation}
\subsection{Matching conditions and sharp interface limit}
We obtain the following matching conditions for the phase field function at $\bm{x} = \bm{\zeta}\left(\bm{s}\right)$: 
\begin{subequations}
    \begin{align}
        \lim_{z\to\pm\infty} \varphi_{0}^{\delta} \left( z, \bm{s} \right) 
        & = \varphi_{0} \left( \bm{x} \pm \right), 
        \label{lim:psi}
        \\ 
        \lim_{z\to \pm\infty} \partial_{z} \varphi_{0}^{\delta} \left( z, \bm{s} \right) & = 0, 
        \label{lim:dpsi0}
        \\ 
        \lim_{z\to\pm\infty} \partial_{z}\varphi_{1}^{\delta}\left( z, \bm{s} \right) 
        & = \nabla \varphi_{0} \left( \bm{x} \pm \right) \cdot \bm{n}, 
        \label{lim:dpsi1}
    \end{align}
\end{subequations}
where $\varphi_0 \left(\bm{x}\pm\right)$ denotes the  $\lim\limits_{\epsilon \searrow 0}\varphi_{0}\left( \bm{x} \pm \epsilon \bm{n} \right)$.  Similar relations hold for the other functions like $\bm{u}, \nu, p$. 

We insert asymptotic expansions \eqref{exp:in} into system \eqref{eqn:main_system} and 
ask that each coefficient of power in $\delta$ vanishes. 
The leading order of equation \eqref{eqn:mu} is $O(1/\delta)$, which gives 
\begin{equation}\label{eqn:psi0}
    0 = \partial_{zz}\varphi_{0}^{\delta} - W^{\prime}\left(\varphi_{0}^{\delta}\right). 
\end{equation}
Thus
\begin{equation}
    0 = \partial_{z}\left(\frac{1}{2}(\partial_{z}\varphi_{0}^{\delta})^2 - W\left(\varphi_{0}^{\delta}\right)\right),\quad \frac{1}{2}(\partial_{z}\varphi_{0}^{\delta})^2 - W\left(\varphi_{0}^{\delta}\right)=C_0, 
\end{equation}
where $C_0$ is a constant in $z$.

From (\ref{lim:psi}), we obtain 
\begin{equation}\label{eqn:limpsi0}
    \varphi_{0}^{\delta} \to \pm 1 \quad \mbox{for} \quad z \to \pm \infty. 
\end{equation}
The above matching condition and (\ref{lim:dpsi0}) yield that $C_0=0$. 
We now assume additionally that 
\begin{equation}\label{eqn:psi}
    \varphi_{0}^{\delta} \left( 0 \right) = 0. 
\end{equation}
By the above conditions \eqref{eqn:psi0}-\eqref{eqn:psi}, $\varphi^{\delta}_{0}$ has a unique solution
\[ \varphi^{\delta}_{0}(z)=\tanh\left(\frac{z}{\sqrt{2}}\right),\]
which does not depend on $t$ and $\bm{s}$. 
The leading order of equation \eqref{eqn:nablau} yields 
\begin{equation}\label{eqn:d_zun}
    \partial_{z} \bm{u}_{0}^{\delta} \cdot \bm{n} 
    = \partial_{z} \left( \bm{u}_{0}^{\delta} \cdot \bm{n} \right) = 0. 
\end{equation}
The matching condition implies that $\left( \bm{u}_{0}^{\delta} \cdot \bm{n} \right) \left( z \right)$ is bounded. Hence 
\begin{equation}
    \left( \bm{u}_{0}^{+} \cdot \bm{n} \right) \left( \bm{x}+ \right) 
    = \lim_{z\to\infty} \left( \bm{u}_{0}^{\delta} \cdot \bm{n} \right) \left( z \right) 
    = \lim_{z\to-\infty} \left( \bm{u}_{0}^{\delta} \cdot \bm{n} \right) \left( z \right) 
    = \left( \bm{u}_{0}^{-} \cdot \bm{n} \right) \left( \bm{x}- \right). 
\end{equation}
Thus 
\begin{equation}
    \left[\bm{u}_{0}\cdot\bm{n}\right]_{-}^{+} = 0, 
\end{equation}
where $\left[f\right]_{-}^{+}\left(\bm{x}\right) = f^{+}\left(\bm{x}\right)-f^{-}\left(\bm{x}\right)$ 
denotes the jump of a quantity at the interface. 
Next, we derive the sharp interface limit for the momentum equation. 
Application of equation \eqref{eqn:inner_grad} for each component yields 
\begin{subequations}
    \begin{align}
        \nabla_{\bm{x}} \bm{u} 
        & = \frac{1}{\delta} \partial_{z} \bm{u} \otimes \bm{n} 
        + \nabla_{\Gamma_{\delta z}} \bm{u}, \label{Eq1}
        \\ 
        \nabla_{\bm{x}} \bm{u} + \nabla_{\bm{x}} \bm{u}^{T} 
        & = \frac{1}{\delta} \left( \partial_{z}\bm{u}\otimes\bm{n} 
        + \bm{n}\otimes\partial_{z}\bm{u} \right) 
        + \left( \nabla_{\Gamma_{\delta z}} \bm{u} 
        + \left( \nabla_{\Gamma_{\delta z}} \bm{u} \right)^{T} \right). \label{Eq2}
    \end{align}
\end{subequations}
We define the $\mathcal{E}\left(\bm{C}\right)=\frac{1}{2}\left(\bm{C}+\bm{C}^{T}\right)$ 
for a quadratic matrix $\bm{C}$. Hence 
\begin{equation}
    \begin{aligned}
    \nabla_{\bm{x}} \cdot 
    \left(\mu_{s}\left(\varphi\right)\left(\nabla_{\bm{x}}\bm{u}
    +\nabla_{\bm{x}}\bm{u}^{T}\right)\right)
    = & \frac{1}{\delta^{2}}\partial_{z}\left( 2 \mu_{s}\left(\varphi^{\delta}\right)\mathcal{E}\left(\partial_{z}\bm{u}^{\delta}\otimes\bm{n}\right)\right)\bm{n}+\frac{1}{\delta}\partial_{z}\left(2\mu_{s}\left(\varphi^{\delta}\right)\mathcal{E}\left(\nabla_{\Gamma_{\delta z}}\bm{u}^{\delta}\right)\right)\bm{n}
    \\
    & +\frac{1}{\delta}\nabla_{\Gamma_{\delta z}}\cdot \left(2\mu_{s}\left(\varphi^{\delta}\right)\mathcal{E}\left(\partial_{z}\bm{u}^{\delta}\otimes\bm{n}\right)\right)+\nabla_{\Gamma_{\delta z}}\cdot\left(2\mu_{s}\left(\varphi^{\delta}\right)\mathcal{E}\left(\nabla_{\Gamma_{\delta z}}\bm{u}^{\delta}\right)\right)
    \\ 
    = & \frac{1}{\delta^{2}}\partial_{z}\left( 2 \mu_{s}\left(\varphi^{\delta}\right)\mathcal{E}\left(\partial_{z}\bm{u}^{\delta}\otimes\bm{n}\right)\bm{n}\right)+\frac{1}{\delta}\partial_{z}\left(2\mu_{s}\left(\varphi^{\delta}\right)\mathcal{E}\left(\nabla_{\Gamma_{\delta z}}\bm{u}^{\delta}\right)\bm{n}\right)
    \\
    & +\frac{1}{\delta}\nabla_{\Gamma_{\delta z}}\cdot \left(2\mu_{s}\left(\varphi^{\delta}\right)\mathcal{E}\left(\partial_{z}\bm{u}^{\delta}\otimes\bm{n}\right)\right)+\nabla_{\Gamma_{\delta z}}\cdot\left(2\mu_{s}\left(\varphi^{\delta}\right)\mathcal{E}\left(\nabla_{\Gamma_{\delta z}}\bm{u}^{\delta}\right)\right), 
    \end{aligned}
\end{equation}
where we used the fact $\partial_{z}\bm{n}=0$. 
We conclude from equation \eqref{eqn:d_zun} 
\begin{equation}
    \left(\bm{n}\otimes\partial_{z}\bm{u}_{0}^{\delta}\right)\bm{n}
    =\left(\partial_{z}\bm{u}_{0}^{\delta}\cdot\bm{n}\right)\bm{n}=\bm{0}.
\end{equation}
Since $\varphi_{0}^{\delta}$ is independent of $t$ and $s$, equation \eqref{eqn:inner_grad} implies 
\begin{equation}
    \nabla \varphi^{\delta} \otimes \nabla \varphi^{\delta} 
     = \frac{1}{\delta^{2}} \left( \partial_{z} \varphi_{0}^{\delta}\right)^{2} \bm{n} \otimes \bm{n} 
     + \frac{2}{\delta} \partial_{z} \varphi_{1}^{\delta} \partial_{z} \varphi_{0}^{\delta} \left( \bm{n} \otimes \bm{n} \right) + \mbox{h.o.t.}. 
\end{equation}
Because of $\left( \nabla_{\Gamma} \bm{n} \right) \bm{n} = \bm{0}$, we conclude 
\begin{equation}
   -\nabla\cdot\bm{\tau}^{\delta}_{\varphi}=\frac{\delta}{Ca} \nabla \cdot \left( \nabla \varphi^{\delta} \otimes \nabla \varphi^{\delta} \right) 
    = \frac{1}{\delta^{2}Ca} \partial_{z} \left(\left( \partial_{z} \varphi_{0}^{\delta} \right)^{2}\right) \bm{n} 
    + \frac{1}{\delta Ca}  \left( \partial_{z} \varphi_{0}^{\delta} \right)^{2} \left( \nabla_{\Gamma} \cdot \bm{n} \right) \bm{n} 
    + \frac{2}{\delta Ca} \partial_{z} \left( \partial_{z} \varphi_{1}^{\delta} \partial_{z} \varphi_{0}^{\delta} \right) \bm{n} + \mbox{h.o.t.}. 
\end{equation}
Since the leading order of the chemical potential does not depend on $z$, 
the term $-\nabla \cdot \left( \bm{u} \otimes M\left( \varphi \right) \nabla \nu \right)$ 
gives no contribution to the order $1/\delta^{2}$~\textcolor{blue}{as  $M\left(\varphi\right)=O(\delta)$}. 
Therefore the order $1/\delta^{2}$ terms from the momentum equation \eqref{eqn:u} yields 
\begin{equation}\label{inner_exp:surface_tension}
     \partial_{z}p^{\delta}_{-1}\bm{n}
    = \partial_{z} \left( \mu_{s} \left( \varphi_{0}^{\delta} \right) \partial_{z} \bm{u}_{0}^{\delta} \right)-\frac{1}{Ca} \partial_{z} \left(\left( \partial_{z} \varphi_{0}^{\delta} \right)^{2}\right)\bm{n}. 
\end{equation}
Taking the inner product of equation \eqref{inner_exp:surface_tension} with $\bm{n}$, 
using $\partial_{z} \bm{n} = 0$ and $\partial_{z} \bm{u}_{0}^{\delta} \cdot \bm{n} = 0$ yields 
\begin{equation}\label{inner_exp:surface_tension2}
   \partial_{z}p^{\delta}_{-1}+\frac{1}{Ca} \partial_{z} \left( \left( \partial_{z} \varphi_{0}^{\delta} \right)^{2} \right) = 0. 
\end{equation}
This gives 
\[ p^{\delta}_{-1}=c-\frac{1}{Ca}\left(\partial_{z} \varphi_{0}^{\delta}\right)^2 \]
where $c$ is a constant in $z$.  The matching condition  $\lim\limits_{z\to \pm\infty}p^{\delta}_{-1}=p_{-1}(\bm{x}\pm)$ gives $c=c^{\pm}$ as $\lim\limits_{z\to \pm\infty}\partial_{z} \varphi_{0}^{\delta}=0.$ As the pressure is not unique up to a constant, therefore, we can choose the  pressures $p$ in the outer regions  are zero at the order  $1/\delta$, i.e. $c^{\pm}=c=0$. In  this case, \[p^{\pm}_{-1}=0,\quad p^{\delta}_{-1}=-\frac{1}{Ca}\left(\partial_{z} \varphi_{0}^{\delta}\right)^2.\]

Equations \eqref{inner_exp:surface_tension} and \eqref{inner_exp:surface_tension2} imply
\begin{equation}\label{inner_exp:viscosity}
    \partial_{z} \left( \mu_{s} \left( \varphi_{0}^{\delta} \right) \partial_{z} \bm{u}_{0}^{\delta} \right) = 0. 
\end{equation}
The matching conditions yield that $\bm{u}_{0}^{\delta}\left(\bm{z}\right)$ is bounded. 
Thus equation \eqref{inner_exp:viscosity} interpreted as an PDE in $z$ has only solutions 
$\bm{u}_{0}^{\delta}$ which are constant in $z$. Again matching yields 
\begin{equation}
    \left[ \bm{u}_{0} \right]_{-}^{+} = 0. 
\end{equation}
Thus we recover the first part of the sharp interface jump condition in equation \eqref{eqn:sharp_condition}. 

Now we analyze the momentum equation to the next order. 
The term $\nabla \cdot \left( \mu_{s} \left( \varphi \right)\left( \nabla \bm{u} + \nabla \bm{u}^{T}\right) \right) $ gives to the order $1/\delta$, 
\begin{equation}
    \partial_{z} \left( \mu_{s} \left( \varphi_{0}^{\delta} \right) \mathcal{E} \left( \partial_{z} \bm{u}_{1}^{\delta} \otimes \bm{n} \right) \bm{n} \right) 
    + \partial_{z} \left( \mu_{s} \left( \varphi_{0}^{\delta} \right) \mathcal{E} \left( \nabla_{\Gamma} \bm{u}_{0}^{\delta} \right) \bm{n} \right). 
\end{equation}
From \eqref{Eq1}, by the matching conditions, we require that
\begin{equation}\label{inner_exp:matching_u}
\partial_{z}\bm{u}^{\delta}_{1}\otimes\bm{n} + \nabla_{\Gamma} \bm{u}^{\delta}_{0} \to \nabla_{\bm{x}} \bm{u}_0 \quad\mbox{for}\quad z\to \pm \infty. 
\end{equation}
From \eqref{eqn:phi}, one has \[\left(\bm{u}^{\delta}_0\cdot \bm{n}-V\right)\partial_z\varphi^{\delta}_0=\frac{M_0}{Pe}\partial_{zz}\nu^{\delta}_0.\]
As $\bm{u}^{\delta}_0$
is independent of $z$, integrating with respect to $z$ and using the matching condition for $\varphi^{\delta}_0, \bm{u}^{\delta}_0$ and $\partial_z\nu^{\delta}_0$,  one has 
\[ [\bm{u}_0\cdot\bm{n}\varphi_0]^{+}_{-}-V[\varphi_0]^{+}_{-}=0.\]
This gives 
\[ \bm{u}^{+}_0\cdot\bm{n}=\bm{u}^{-}_0\cdot\bm{n}=V,\]
as $[\varphi_0]^{+}_{-}=2.$
By matching conditions, this yields 
\[ \bm{u}^{\delta}_0\cdot\bm{n}-V=0.\]
Therefore,
\[\partial_{zz}\nu^{\delta}_0=0.\]
This gives $\nu^{\delta}_0$ is a constant in $z$ by using the matching condition.

Moreover, the term $-\nabla\cdot\left(\bm{u}\otimes M\left(\varphi\right)\nabla\nu\right)$ 
gives no contribution to order $1/\delta$ as $M\left(\varphi\right)$ is a constant and $M\left(\varphi\right)=O(\delta)$. 
Thus we obtain from the momentum equation at order $1/\delta$: 
\begin{equation}
    \begin{aligned}
        & \mathrm{Re} \left( -\partial_{z}\left(\rho\left(\varphi_{0}^{\delta}\right)\bm{u}_{0}^{\delta}\right) V + \partial_{z} \left( \rho \left( \varphi_{0}^{\delta} \right) \left( \bm{u}_{0}^{\delta} \otimes \bm{u}_{0}^{\delta} \right)\right) \bm{n} \right) 
        \\ 
        & - 2 \partial_{z}\left(\mu_{s} \left( \varphi_{0}^{\delta} \right) \mathcal{E} \left(\partial_{z}\bm{u}_{1}^{\delta} \otimes \bm{n} \right) \bm{n} \right) 
        - 2\partial_{z} \left( \mu_{s} \left( \varphi_{0}^{\delta} \right) \mathcal{E} \left( \nabla_{\Gamma} \bm{u}_{0}^{\delta} \right) \bm{n} \right) 
        \\ 
        & - \frac{1}{\mathrm{Ca}} \left( \partial_{z} \varphi_{0}^{\delta} \right)^{2} \kappa \bm{n} 
        + \frac{1}{\mathrm{Ca}} \partial_{z} \left( 2 \partial_{z} \varphi_{1}^{\delta} \partial_{z} \varphi_{0}^{\delta} \right) \bm{n} 
        -\partial_{z} \left( \bm{\tau}_{p,0}^{\delta} \right) \bm{n} + \partial_{z} p_{0}^{\delta} \bm{n} = 0, 
    \end{aligned}
\end{equation}
because of $\kappa = - \nabla\cdot \bm{n}$. 
Integration with respect to $z$ yields after matching and the use of equation \eqref{inner_exp:matching_u} 
\begin{equation}
    \begin{aligned}
    & \mathrm{Re} \left( - \left[ \rho \bm{u}_{0} \right]_{-}^{+}V
    + \left[ \rho \bm{u}_{0}\bm{u}_{0} \cdot \bm{n} \right]_{-}^{+}  \right) 
    - 2 \left[ \mu_{s} \mathcal{E}\left(\nabla_{\bm{x}}\bm{u}_{0} \right) \right]_{-}^{+}\bm{n} 
  - \frac{1}{\mathrm{Ca}} \int_{-\infty}^{\infty} \left( \partial_{z}\varphi_{0}^{\delta}\right)^{2} \mathrm{d} \bm{z} \kappa \bm{n} 
    -\left[ \bm{\tau}_{p,0} \right]_{-}^{+} 
    + \left[ p_{0} \right]_{-}^{+} \bm{n} = 0, 
    \end{aligned}
\end{equation}
where $\int_{-\infty}^{\infty} \left( \partial_{z} \varphi_{0}^{\delta}\right)^{2} \mathrm{d}\bm{z} = \frac{2\sqrt{2}}{3}$ for the choice of $W$,  $\rho^{\pm}=\rho(\varphi^{\pm}_0)=\rho(\pm 1)$ and $\mu_s^{\pm}=\mu_s(\varphi^{\pm}_0)=\mu_s(\pm 1)$ are used. 
Since $\bm{u}^+_{0} \cdot \bm{n}=\bm{u}^-_{0} \cdot \bm{n} = V, [\bm{u}_0]^+_{-}=0$, we conclude 
\begin{equation}
    -2\left[ \mu_{s} D\bm{u}_{0} \right]_{-}^{+} \bm{n} 
    + \left[ p_{0} \right]_{-}^{+} \bm{n} 
    -\left[ \bm{\tau}_{p,0} \right]_{-}^{+} \bm{n} 
    = \frac{2\sqrt{2}}{3} \frac{1}{\mathrm{Ca}} \kappa \bm{n}. 
\end{equation} 
Therefore the leading-order equaton recovers the dimensionless version of the sharp interface jump condition in equation \eqref{eqn:sharp_condition}. 
\section[Time discrete scheme and energy stability]{Time discrete scheme and energy stability}\label{sec4}
Before conducting the numerical experiments, we construct a numerical scheme
for the proposed FENE-P--CH model~\eqref{eqn:main_system}.
Our aim is to preserve the energy dissipation law in Theorem~\ref{thm:energy}
after time discretization. In this section, we discretize only the time
variable and retain the continuous spatial operators. 
We first rewrite the system in a form suitable for the energy estimate, 
then introduce the time discrete scheme and prove its unconditional energy stability. 

For convenience, we use the following notations:
\begin{equation}
 W(\varphi)=\frac14(\varphi^2-1)^2,\qquad
 m=\frac{M(\varphi)}{\Pe}=\frac{M_0\delta}{\Pe},\qquad
 \rho(\varphi)=\frac{\rho^++\rho^-}{2}
                 +\frac{\rho^+-\rho^-}{2}\varphi.
 \label{cs:eq:rho}
\end{equation}
We write $\rho'=(\rho^+-\rho^-)/2$ and define
\begin{equation}
 q(\varphi)=\frac{\beta\mu_p(\varphi)}{2\Wi\lambda(\varphi)},\qquad
 r(\varphi)=\frac{1}{\Wi\lambda(\varphi)}.
 \label{cs:eq:qr}
\end{equation}
Throughout this section, $\Rey,\Ca,\Wi,\delta,m>0$, $\beta\geq0$,
$\mu_s(\varphi)>0$, $\mu_p(\varphi)\geq0$, and $\lambda(\varphi)>0$.
The material functions are continuously differentiable on the attained
phase-field range. In particular, $q\geq0$ and $r>0$.

Let $\Omega\subset\mathbb R^d$, where $d=2$ or $3$, and let $L^2>d$.
The admissible conformation tensors belong to
\begin{equation}
 \mathcal K=\{\bm\sigma=\bm\sigma^T\succ0:
                  \tr\bm\sigma<L^2\},
 \label{cs:eq:admissible}
\end{equation}
where $\mathcal K$ is set of symmetric positive definite matrices whose trace is less than $L^2$.

For $\bm\sigma\in\mathcal K$, set
\begin{align}
 g(\bm\sigma)&=\left(1-\frac{\tr\bm\sigma}{L^2}\right)^{-1},
 & f(\bm\sigma)&=g(\bm\sigma)\bm\sigma-\bm I,
 \label{cs:eq:fs}\\
 H(\bm\sigma)&=-L^2\ln\left(1-\frac{\tr\bm\sigma}{L^2}\right)-\tr(\ln\bm\sigma)-\tr(\bm I),
 & \bm B(\bm\sigma)&=g(\bm\sigma)\bm I-\bm\sigma^{-1}.
 \label{cs:eq:entropy}
\end{align}
With  $q$ defined in~\eqref{cs:eq:qr}, the elastic energy density is
$q(\varphi)H(\bm\sigma)$, 
and the polymer stress is $\tau_p=2q(\varphi)f(\bm\sigma)$.

We denote the $L^2(\Omega)$ inner product by $(\cdot,\cdot)$ and the
corresponding norm by $\|\cdot\|$. For vector and tensor fields, the
integrands use the Euclidean dot product and Frobenius contraction,
respectively. We write $\bm A=\nabla\bm u$, with
$A_{ij}=\partial_j u_i$, so that the stretching term is
$\bm A\bm\sigma+\bm\sigma\bm A^T$ and
$(\nabla\cdot\bm T)_i=\sum\limits_j\partial_j T_{ij}$.
We impose periodic boundary conditions, or the folllowing homogeneous conditions
\begin{equation}
 \bm u|_{\partial\Omega}=0,\qquad
 \partial_{\bm n}\varphi|_{\partial\Omega}=0,\qquad
 \partial_{\bm n}\nu|_{\partial\Omega}=0.
 \label{cs:eq:boundary}
\end{equation}
These conditions eliminate the boundary terms in the integrations by parts
below. No conformation diffusion is introduced.

Define the total mass flux and the transport operator by
\begin{equation}
 \bm a=\rho\bm u-\rho'm\nabla\nu,\qquad
 \mathcal B(\bm a,\bm u)
 =\bm a\cdot\nabla\bm u+\frac12(\nabla\cdot\bm a)\bm u.
 \label{cs:eq:skewdef}
\end{equation}
The phase-field equation and incompressibility give
$\partial_t\rho+\nabla\cdot\bm a=0$. Consequently,
\begin{equation}
 \rho\partial_t\bm u+\bm a\cdot\nabla\bm u
 =\rho\partial_t\bm u+\frac12(\partial_t\rho)\bm u
      +\mathcal B(\bm a,\bm u).
 \label{cs:eq:continuous-inertia}
\end{equation}
Moreover, integration by parts yields
$(\mathcal B(\bm a,\bm u),\bm u)=0$.
The capillary and elastic terms can be written as
$-\varphi\nabla\nu/\Ca+q\nabla H(\bm\sigma)$ after absorbing the
gradient potentials into the pressure. Indeed, with
$e_{\rm mix}=\delta|\nabla\varphi|^2/2+W(\varphi)/\delta$,
\begin{equation}
 \frac1\Ca\nabla\cdot(-\delta\nabla\varphi\otimes\nabla\varphi)
 =-\frac{\varphi\nabla\nu}{\Ca}+q\nabla H(\bm\sigma)
   +\nabla\left(\frac{\varphi\nu-e_{\rm mix}}{\Ca}
                         -qH(\bm\sigma)\right).
 \label{cs:eq:pressure-rewrite}
\end{equation}
Hereafter, $P=p-\frac{\varphi\nu-e_{\rm mix}}{\Ca}
                         +qH(\bm\sigma)$ denotes the modified pressure.

Now we rewrite the governing system  \eqref{eqn:main_system}  as
\begin{subequations}\label{cs:eq:equivalent}
\begin{align}
 \partial_t\varphi+\bm u\cdot\nabla\varphi
     &=m\Delta\nu,\\
 \nu&=-\delta\Delta\varphi+\frac{W'(\varphi)}{\delta}
                           +\Ca q'(\varphi)H(\bm\sigma),\\
 \Rey\left(\rho\partial_t\bm u+\frac12(\partial_t\rho)\bm u
                      +\mathcal B(\bm a,\bm u)\right)+\nabla P
     &=\nabla\cdot\left(\mu_s(\varphi)(\bm A+\bm A^T)+2q(\varphi)f(\bm\sigma)\right)
       \nonumber\\
     &\quad-\frac{\varphi\nabla\nu}{\Ca}
                +q(\varphi)\nabla H(\bm\sigma),\\
 \nabla\cdot\bm u&=0,\\
 \partial_t\bm\sigma+\bm u\cdot\nabla\bm\sigma
                 -\bm A\bm\sigma-\bm\sigma\bm A^T
     &=-r(\varphi) f(\bm\sigma). 
\end{align}
\end{subequations}

\begin{remark}
From the above asymptotic analysis, the original pressure $p$ in the inner layer has the order $O(1/\delta)$. However, the modified pressure $P$ in the inner layer at the order $1/\delta$ is zero since
\[p^{\delta}_{-1}+\frac{1}{Ca}\left(\frac{1}{2}(\partial_z \phi^{\delta}_0)^2+W(\phi^{\delta}_0)\right)=p^{\delta}_{-1}+\frac{(\partial_z \phi^{\delta}_0)^2}{Ca}=0.\]
\end{remark}

\renewcommand{\qedsymbol}{}
\subsection{Semi-discretization in time}\label{cs:sec:scheme}\label{sec5}
Let $t^n$ denote the time levels and let
$\tau_n=t^{n+1}-t^n>0$. The superscripts $n$ and $n+1$ denote the known
and new states, respectively. For the material functions, we use the notation
$\rho^k=\rho(\varphi^k)$, $q^k=q(\varphi^k)$,
$r^k=r(\varphi^k)$,  and $\mu_s^k=\mu_s(\varphi^k)$ ($k=n,n+1$).
Define the divided difference
\begin{equation}
 q(\varphi^{n+1},\varphi^n)=
 \begin{cases}
 \dfrac{q^{n+1}-q^n}{\varphi^{n+1}-\varphi^n},
                                  &\varphi^{n+1}\ne\varphi^n,\\[4pt]
 q'(\varphi^n),&\varphi^{n+1}=\varphi^n.
 \end{cases}
 \label{cs:eq:secant}
\end{equation}

Suppose that $(\varphi^n,\bm u^n,\bm\sigma^n)$ are known.
We seek $(\varphi^{n+1},\nu^{n+1},\bm u^{n+1},P^{n+1},\bm\sigma^{n+1})$ 
satisfying the following coupled equations:
\begin{subequations}\label{cs:eq:ch-sigma}
\begin{align}
 \frac{\varphi^{n+1}-\varphi^n}{\tau_n}
       +\bm u^{n+1}\cdot\nabla\varphi^{n+1}
     &=m\Delta\nu^{n+1},\label{cs:eq:phase}\\
 \nu^{n+1}& =-\delta\Delta\varphi^{n+1}
       +\frac{(\varphi^{n+1})^3-\varphi^n}{\delta}
       +\Ca q(\varphi^{n+1},\varphi^n)H(\bm\sigma^n),
       \label{cs:eq:chemical}
\end{align}
\begin{equation}
\begin{aligned}
 &\frac{\bm\sigma^{n+1}-\bm\sigma^n}{\tau_n}
       +\bm u^{n+1}\cdot\nabla\bm\sigma^{n+1}
       -\bm A^{n+1}\bm\sigma^{n+1}
       -\bm\sigma^{n+1}(\bm A^{n+1})^T
     =-r^{n+1}f(\bm\sigma^{n+1}).
\end{aligned}
\label{cs:eq:sigma}
\end{equation}
\end{subequations}
where $\bm A^{n+1}=\nabla\bm u^{n+1}$. Set
\begin{equation}
 \bm a^{n+1}=\rho^{n+1}\bm u^{n+1}-\rho'm\nabla\nu^{n+1}.
 \label{cs:eq:flux}
\end{equation}
The momentum equation and incompressibility constraint are
\begin{subequations}\label{cs:eq:momentum-system}
\begin{align}
 &\Rey\left(\rho^n\frac{\bm u^{n+1}-\bm u^n}{\tau_n}
       +\frac{\rho^{n+1}-\rho^n}{2\tau_n}\bm u^{n+1}
       +\mathcal B(\bm a^{n+1},\bm u^{n+1})\right)+\nabla P^{n+1}
       \nonumber\\
 &\quad=\nabla\cdot\left(\mu_s^{n+1}
                (\bm A^{n+1}+(\bm A^{n+1})^T)
                +2q^{n+1}f(\bm\sigma^{n+1})\right)
       -\frac{\varphi^{n+1}\nabla\nu^{n+1}}{\Ca}
            +q^{n+1}\nabla H(\bm\sigma^{n+1}),\label{cs:eq:momentum}\\
 &\nabla\cdot\bm u^{n+1}=0,\quad
                  \int_\Omega P^{n+1}\,\mathrm d\bm x=0.
                  \label{cs:eq:incompressible}
\end{align}
\end{subequations}
The periodic boundary conditions or boundary conditions \eqref{cs:eq:boundary} are imposed at each time
level. 
The quartic part of $W$ is treated implicitly and its concave quadratic part explicitly, 
following the convex-splitting principle~\cite{Eyre1998}. 
The divided difference accounts exactly for changes in the phase-dependent elastic coefficient. 

Integrating~\eqref{cs:eq:phase} over $\Omega$ and using the fact $\rho^k=\rho(\varphi^k)$ ($k=n, n+1$), one can get the  following mass conservation:
\begin{equation}
 \int_\Omega\varphi^{n+1}\,\mathrm d\bm x
   =\int_\Omega\varphi^n\,\mathrm d\bm x,\quad
 \int_\Omega\rho^{n+1}\,\mathrm d\bm x
   =\int_\Omega\rho^n\,\mathrm d\bm x.
 \label{cs:eq:mass}
\end{equation}
The affine density law also yields
\begin{equation}
 \frac{\rho^{n+1}-\rho^n}{\tau_n}
                       +\nabla\cdot\bm a^{n+1}=0.
 \label{cs:eq:mass-local}
\end{equation}

\subsection{Unconditional energy stability}\label{cs:sec:theorem}
Before proving the energy estimate, we introduce the following identities.
\begin{lemma}\label{cs:lem:algebra}
For any real numbers $x,y$ and vectors $\bm v,\bm w$,
\begin{align}
 2\bm v\cdot(\bm v-\bm w)
   &=|\bm v|^2-|\bm w|^2+|\bm v-\bm w|^2,
       \label{cs:eq:polarization}\\
 (x^3-y)(x-y)&=W(x)-W(y)+\mathcal R_W(x,y),
       \label{cs:eq:quartic}\\
 \mathcal R_W(x,y)&=\frac{(x-y)^2}{4}(3x^2+2xy+y^2+2)
                         \geq\frac12(x-y)^2.\nonumber
\end{align}
For positive numbers $\rho^{n+1},\rho^n$, one also has
\begin{equation}
 \left(\rho^n(\bm u^{n+1}-\bm u^n)
        +\frac{\rho^{n+1}-\rho^n}{2}\bm u^{n+1}\right)
                     \cdot\bm u^{n+1}
 =\frac12\left(\rho^{n+1}|\bm u^{n+1}|^2
               -\rho^n|\bm u^n|^2
               +\rho^n|\bm u^{n+1}-\bm u^n|^2\right).
 \label{cs:eq:kinetic}
\end{equation}
\end{lemma}
\begin{proof}
The identities follow by expanding both sides. The lower bound for
$\mathcal R_W$ follows from $3x^2+2xy+y^2=(x+y)^2+2x^2\geq0$.
\end{proof}

\begin{lemma}\label{cs:lem:entropy}
For $\bm\sigma,\bm\sigma_0\in\mathcal K$, the function $H$ satisfies
$H'(\bm\sigma)=\bm B(\bm\sigma)$ and
\begin{equation}
 \mathcal R_H(\bm\sigma,\bm\sigma_0)
 :=\bm B(\bm\sigma):(\bm\sigma-\bm\sigma_0)
                  -H(\bm\sigma)+H(\bm\sigma_0)\geq0.
 \label{cs:eq:RH}
\end{equation}
Furthermore, for any real matrix $\bm A$,
\begin{align}
 \bm B(\bm\sigma):(\bm A\bm\sigma+\bm\sigma\bm A^T)
       &=2f(\bm\sigma):\bm A,\label{cs:eq:stretch}\\
 \bm B(\bm\sigma):f(\bm\sigma)
       &=\tr\big(f(\bm\sigma)^2\bm\sigma^{-1}\big)\geq0.
       \label{cs:eq:relax}
\end{align}
\end{lemma}
\begin{proof}
For every nonzero symmetric matrix $\bm Z$, direct differentiation gives
\begin{equation}
 H''(\bm\sigma)(\bm Z,\bm Z)
  =\frac{g(\bm\sigma)^2}{L^2}(\tr\bm Z)^2
       +\tr(\bm\sigma^{-1}\bm Z\bm\sigma^{-1}\bm Z)>0.
 \label{cs:eq:hessian}
\end{equation}
The second term is $\|\bm\sigma^{-1/2}\bm Z\bm\sigma^{-1/2}\|_F^2$.
Thus $H$ is strictly convex on the convex set $\mathcal K$, which proves
\eqref{cs:eq:RH}. The relation
$\bm B\bm\sigma=\bm\sigma\bm B=f(\bm\sigma)$ proves
\eqref{cs:eq:stretch} by cyclic invariance of the trace.
If $\ell_1,\ldots,\ell_d$ are the positive eigenvalues of $\bm\sigma$,
then the last expression in~\eqref{cs:eq:relax} equals
$\sum\limits_{j=1}^d(f(\ell_j)-1)^2/\ell_j\geq0$.
\end{proof}

Define the total energy at time level $n$ by
\begin{equation}
 E_{\rm total}^n=E^n_{kin}+E^n_{mix}+E^n_{el},
\label{cs:eq:energy}
\end{equation}
where $E^n_{kin}, E^n_{mix}, E^n_{el}$ are the kinetic, mixing, and elastic energy of the original model,
evaluated at the  time-discrete level $n$, given respectively by 
\[E^n_{kin}=\frac{\Rey}{2}\|\sqrt{\rho^n}\bm u^n\|^2,  \quad E^n_{mix}=\frac1\Ca\int_\Omega \left(\frac\delta2|\nabla\varphi^n|^2+\frac{W(\varphi^n)}\delta\right)\mathrm d\bm x, \quad  E^n_{el}=\int_\Omega q^nH(\bm\sigma^n) \mathrm d\bm x.\]
The discrete energy dissipation rate and the work done by  external force are given respectively by
\begin{align}
 \mathcal D^{n+1}={}&\frac m\Ca\|\nabla\nu^{n+1}\|^2
   +\frac12\|\sqrt{\mu_s^{n+1}}
                   (\bm A^{n+1}+(\bm A^{n+1})^T)\|^2
                   +\int_\Omega q^{n+1}r^{n+1}
      \tr\big(f(\bm\sigma^{n+1})^2(\bm\sigma^{n+1})^{-1}\big)
           \mathrm d\bm x,\label{cs:eq:dissipation}. 
\end{align}

\begin{theorem}[Unconditional energy stability]\label{cs:thm:energy}
Suppose that the coefficient assumptions of Section~\ref{sec4}
hold and that the time-discrete equations
\eqref{cs:eq:ch-sigma}--\eqref{cs:eq:momentum-system} admit a sufficiently
regular solution satisfying $\rho^n,\rho^{n+1}>0$ and
$\bm\sigma^n,\bm\sigma^{n+1}\in\mathcal K$ throughout $\Omega$.
Under the stated periodic or homogeneous boundary conditions~(\ref{cs:eq:boundary}), the scheme
satisfies, for every $\tau_n>0$,
\begin{equation}
 E_{\rm total}^{n+1}-E_{\rm total}^n
     \leq-\tau_n\mathcal D^{n+1}.
 \label{cs:eq:energy-inequality}
\end{equation}
\end{theorem}

\begin{proof}
We first take the $L^2$ inner product of~\eqref{cs:eq:phase} with
$\tau_n\nu^{n+1}/\Ca$. Integration by parts, the boundary conditions,
and $\nabla\cdot\bm u^{n+1}=0$ give
\begin{equation}
 \frac1\Ca(\varphi^{n+1}-\varphi^n,\nu^{n+1})
  =-\frac{\tau_n m}{\Ca}\|\nabla\nu^{n+1}\|^2
   +\frac{\tau_n}{\Ca}
           (\varphi^{n+1}\nabla\nu^{n+1},\bm u^{n+1}).
\label{cs:eq:phase-test}
\end{equation}
Similarly, taking the $L^2$ inner product of~\eqref{cs:eq:chemical}
with $(\varphi^{n+1}-\varphi^n)/\Ca$ yields
\begin{equation}
\begin{aligned}
 \frac1\Ca(\varphi^{n+1}-\varphi^n,\nu^{n+1})
 ={}&\frac\delta\Ca
     (\nabla\varphi^{n+1},\nabla\varphi^{n+1}-\nabla\varphi^n)
     +\frac1{\Ca\delta}
       ((\varphi^{n+1})^3-\varphi^n,\varphi^{n+1}-\varphi^n)\\
 &+\int_\Omega(q^{n+1}-q^n)H(\bm\sigma^n) {\mathrm d}\bm {x}.
\end{aligned}
\label{cs:eq:chemical-test}
\end{equation}
Here we have used
$q(\varphi^{n+1},\varphi^n)(\varphi^{n+1}-\varphi^n)=q^{n+1}-q^n$, 
including at points where $\varphi^{n+1}=\varphi^n$. Applying Lemma~\ref{cs:lem:algebra}
to~\eqref{cs:eq:chemical-test} and then using~\eqref{cs:eq:phase-test}, we obtain
\begin{equation}
\begin{aligned}
 &E_{\rm mix}^{n+1}-E_{\rm mix}^n
  +\frac\delta{2\Ca}\|\nabla\varphi^{n+1}-\nabla\varphi^n\|^2
  +\frac1{\Ca\delta}\int_\Omega
          \mathcal R_W(\varphi^{n+1},\varphi^n) {\mathrm d}\bm {x}
          +\int_\Omega(q^{n+1}-q^n)H(\bm\sigma^n) {\mathrm d}\bm {x}\\
  &=-\frac{\tau_n m}{\Ca}\|\nabla\nu^{n+1}\|^2
   +\frac{\tau_n}{\Ca}
         (\varphi^{n+1}\nabla\nu^{n+1},\bm u^{n+1}).
\end{aligned}
\label{cs:eq:mixing-estimate}
\end{equation}

Next, we take the tensor $L^2$ inner product of~\eqref{cs:eq:sigma}
with $\tau_n q^{n+1}\bm B(\bm\sigma^{n+1})$.
The chain rule in the spatial variable gives
\begin{equation}
 \bm B(\bm\sigma^{n+1}):
       (\bm u^{n+1}\cdot\nabla\bm\sigma^{n+1})
   =\bm u^{n+1}\cdot\nabla H(\bm\sigma^{n+1}).
 \label{cs:eq:transport-cancel}
\end{equation}
Together with Lemma~\ref{cs:lem:entropy}, this yields
\begin{equation}
\begin{aligned}
 &\int_\Omega(\bm\sigma^{n+1}-\bm\sigma^n):
                      (q^{n+1}\bm B(\bm\sigma^{n+1})){\mathrm d}\bm {x}
   +\tau_n(q^{n+1}\nabla H(\bm\sigma^{n+1}),\bm u^{n+1})
   -2\tau_n\int_{\Omega}(q^{n+1}f(\bm\sigma^{n+1})):\bm A^{n+1}{\mathrm d}\bm {x}\\
   &=-\tau_n\int_\Omega q^{n+1}r^{n+1}\tr\big(f(\bm\sigma^{n+1})^2(\bm\sigma^{n+1})^{-1}\big)
         {\mathrm d}\bm {x}.
\end{aligned}
\label{cs:eq:sigma-test}
\end{equation}

By the definition of $\mathcal R_H$, the first term satisfies
\begin{equation}
\int_\Omega(\bm\sigma^{n+1}-\bm\sigma^n):
                      (q^{n+1}\bm B(\bm\sigma^{n+1})){\mathrm d}\bm {x}
  =E_{\rm el}^{n+1}-E_{\rm el}^n
       -\int_\Omega(q^{n+1}-q^n)H(\bm\sigma^n) {\mathrm d}\bm {x}
       +\int_\Omega q^{n+1}
              \mathcal R_H(\bm\sigma^{n+1},\bm\sigma^n) {\mathrm d}\bm {x}.
\label{cs:eq:elastic-product}
\end{equation}

Finally, we take the $L^2$ inner product of~\eqref{cs:eq:momentum}
with $\tau_n\bm u^{n+1}$. The transport and pressure terms vanish:
\begin{align}
 (\mathcal B(\bm a^{n+1},\bm u^{n+1}),\bm u^{n+1})&=0,
       \label{cs:eq:skewzero}\\
 (\nabla P^{n+1},\bm u^{n+1})
       &=-(P^{n+1},\nabla\cdot\bm u^{n+1})=0.\nonumber
\end{align}
Integration by parts in the full viscous stress gives
\begin{equation}
 \left(\nabla\cdot\left(\mu_s^{n+1}
       (\bm A^{n+1}+(\bm A^{n+1})^T)\right),\bm u^{n+1}\right)
  =-\frac12\|\sqrt{\mu_s^{n+1}}
       (\bm A^{n+1}+(\bm A^{n+1})^T)\|^2.
 \label{cs:eq:viscous}
\end{equation}
The identity~\eqref{cs:eq:kinetic}  gives
\begin{equation}
\begin{aligned}
 &E_{\rm kin}^{n+1}-E_{\rm kin}^n
       +\frac\Rey2\|\sqrt{\rho^n}(\bm u^{n+1}-\bm u^n)\|^2\\
 &\quad=-\frac{\tau_n}{2}\|\sqrt{\mu_s^{n+1}}
              (\bm A^{n+1}+(\bm A^{n+1})^T)\|^2
       -2\tau_n\int_{\Omega}(q^{n+1}f(\bm\sigma^{n+1})):\bm A^{n+1}{\mathrm d}\bm{x}\\
 &\quad-\frac{\tau_n}{\Ca}
         (\varphi^{n+1}\nabla\nu^{n+1},\bm u^{n+1})
       +\tau_n(q^{n+1}\nabla H(\bm\sigma^{n+1}),\bm u^{n+1}).
\end{aligned}
\label{cs:eq:momentum-test}
\end{equation}

Substituting~\eqref{cs:eq:elastic-product} into~\eqref{cs:eq:sigma-test}
and adding the resulting identity to~\eqref{cs:eq:mixing-estimate}
and~\eqref{cs:eq:momentum-test}, we find that the capillary work,
polymer-stress work, conformation-transport work, and variation of the
elastic coefficient cancel. Hence,
\begin{equation}
 E_{\rm total}^{n+1}-E_{\rm total}^n+\mathcal R^{n+1}
           +\tau_n\mathcal D^{n+1}=0,
 \label{cs:eq:exact-energy}
\end{equation}
where
\begin{equation}
\begin{aligned}
 \mathcal R^{n+1}={}&\frac\Rey2
       \|\sqrt{\rho^n}(\bm u^{n+1}-\bm u^n)\|^2
       +\frac\delta{2\Ca}
                  \|\nabla\varphi^{n+1}-\nabla\varphi^n\|^2\\
 &+\frac1{\Ca\delta}\int_\Omega
             \mathcal R_W(\varphi^{n+1},\varphi^n) \mathrm d\bm x
       +\int_\Omega q^{n+1}
             \mathcal R_H(\bm\sigma^{n+1},\bm\sigma^n) \mathrm d\bm x
       \geq0.
\end{aligned}
\label{cs:eq:remainder}
\end{equation}
Dropping these nonnegative terms proves~\eqref{cs:eq:energy-inequality}.
The proof uses no upper bound on the time step or on $W''$, and is valid
when $q=0$ in a Newtonian phase. This completes the proof.
\end{proof}

\begin{remark}
The theorem gives an energy estimate for every admissible, sufficiently
regular solution of the time-discrete equations. It does not assert
the numerical preservation for the positivity of
density and  the positive definiteness of the conformation tensor. For forced flows, the estimate includes
the work of the external force, so the total energy need not decrease.
\end{remark}

\section{Numerical experiments}\label{sec6}

In this section, we conduct a series of two-dimensional numerical experiments to 
examine temporal and spatial convergence, energy dissipation,
and droplet deformation in periodic and boundary-driven shear flows. 
All experiments use unequal bulk densities $\rho^+=1$ and $\rho^-=1.5$, and the domain is assumed to be $\Omega=[-1,1]^2$. 

The time discrete scheme and its energy estimate in Section~\ref{sec4}
use continuous spatial operators. For spatial discretization, we mainly focus on periodic boundary conditions where the Fourier spectral  collocation method is used. In addition, in order to further illustrate the effectiveness of the proposed method, we also investigate the Couette flow, where the domain is periodic in $x$ direction and bounded by two moving walls at $y=\pm1$ with opposite velocities, and a combination of  Fourier and Legendre collocation methods is used.  The temporal convergence, energy, body force shear, and Couette calculations
use $129\times129$ grid points. Only the spatial-convergence study changes this
grid, with its finest resolution also $129\times129$.

\subsection{Temporal and spatial convergence}

For this study, periodic boundary conditions are used. We choose the parameters 
$\mathrm{Re}=\mathrm{Ca}=\mathrm{Wi}=\beta=1$, $L^{2}=20$, $\delta=0.08$, $m=0.002$.
The bulk properties are
$(\rho^{+},\rho^{-})=(1,1.5)$,
$(\mu^{+}_{s},\mu^{-}_{s})=(0.7,1)$,
$(\mu^{+}_{p},\mu^{-}_{p})=(0.5,0.3)$, and
$(\lambda^{+},\lambda^{-})=(0.7,1)$.
All grids sample the same smooth periodic initial data:
\begin{align}
\varphi_0&=0.25+0.12\sin(\pi x)\cos(\pi y)
             +0.04e^{\cos(\pi x)+\cos(\pi y)-2},\nonumber\\
 \bm{u}_0&=\left(u_{x,0}, u_{y,0}\right)^T=\left(0.04\sin(\pi x)\cos(\pi y),-0.04\cos(\pi x)\sin(\pi y)\right)^T,\nonumber\\
 \bm{\sigma}_0&= \begin{pmatrix}
         \sigma_{xx,0} \quad  \sigma_{xy,0} \\
         \sigma_{xy,0} \quad  \sigma_{yy,0}
    \end{pmatrix}= \begin{pmatrix}
         c_{\rm eq}+0.035\cos(\pi x)\cos(\pi y) \quad 0.0175\sin(\pi x)\sin(\pi y) \\
         0.0175\sin(\pi x)\sin(\pi y) \quad c_{\rm eq}-0.035\cos(\pi x)\cos(\pi y)  
    \end{pmatrix}
\label{eq:spectral_initial}
\end{align}
with $c_{\rm eq}=L^{2}/(L^{2}+2)$. 
The equal velocity amplitudes ensure $\nabla\cdot\bm u_0=0$.
For temporal convergence, we fix $N=129$, $t=0.08$, and use
$\tau=0.01,0.005,0.0025,0.00125,0.000625$.
For each variable $U$, we compute
\begin{align*}
 e_\tau(U)=\normh{U_\tau(0.08)-U_{\tau/2}(0.08)},\quad
 p_\tau(U)=\log_2\frac{e_{2\tau}(U)}{e_\tau(U)}.
\end{align*}
Table~\ref{tab:spectral_time} reports the temporal order as $\tau$ decreases.

\begin{table}[htbp]
\centering\small
\setlength{\tabcolsep}{3pt}
\caption{Temporal $L^2$ Cauchy errors and orders at $t=0.08$ on the $129\times129$ spatial grid.}
\label{tab:spectral_time}
\begin{tabular}{lrrrrrrrr}
\toprule
Variable & \multicolumn{2}{c}{$\tau=0.01$} & \multicolumn{2}{c}{$\tau=0.005$} & \multicolumn{2}{c}{$\tau=0.0025$} & \multicolumn{2}{c}{$\tau=0.00125$} \\
\cmidrule(lr){2-3}\cmidrule(lr){4-5}\cmidrule(lr){6-7}\cmidrule(lr){8-9}
 & Error & Order & Error & Order & Error & Order & Error & Order \\
\midrule
$\varphi$ & $1.374\times10^{-5}$ & --- & $7.093\times10^{-6}$ & 0.9541 & $3.605\times10^{-6}$ & 0.9763 & $1.818\times10^{-6}$ & 0.9880 \\
$u_x$ & $5.128\times10^{-4}$ & --- & $2.663\times10^{-4}$ & 0.9456 & $1.357\times10^{-4}$ & 0.9721 & $6.853\times10^{-5}$ & 0.9859 \\
$u_y$ & $5.128\times10^{-4}$ & --- & $2.662\times10^{-4}$ & 0.9456 & $1.357\times10^{-4}$ & 0.9721 & $6.852\times10^{-5}$ & 0.9859 \\
$\sigma_{xx}$ & $2.179\times10^{-4}$ & --- & $1.137\times10^{-4}$ & 0.9381 & $5.813\times10^{-5}$ & 0.9684 & $2.939\times10^{-5}$ & 0.9840 \\
$\sigma_{xy}$ & $6.186\times10^{-6}$ & --- & $3.124\times10^{-6}$ & 0.9856 & $1.570\times10^{-6}$ & 0.9926 & $7.871\times10^{-7}$ & 0.9962 \\
$\sigma_{yy}$ & $2.179\times10^{-4}$ & --- & $1.137\times10^{-4}$ & 0.9381 & $5.813\times10^{-5}$ & 0.9684 & $2.939\times10^{-5}$ & 0.9840 \\
$\nu$ & $2.644\times10^{-3}$ & --- & $1.327\times10^{-3}$ & 0.9946 & $6.647\times10^{-4}$ & 0.9973 & $3.327\times10^{-4}$ & 0.9986 \\
$P$ & $6.252\times10^{-4}$ & --- & $3.136\times10^{-4}$ & 0.9956 & $1.570\times10^{-4}$ & 0.9978 & $7.857\times10^{-5}$ & 0.9989 \\
\bottomrule
\end{tabular}
\end{table}

The independent spatial study fixes $\tau=0.0005$ and $t=0.02$,
and uses $N=9,17,33,65,129$. 
We compare
$e_N(U)=\normh{U_N-\mathcal R_NU_{N_{\rm fine}}}$,
where $\mathcal R_N$ truncates the normalized Fourier coefficients of the fine solution to the coarse mode set. 
Table~\ref{tab:spectral_space} shows decreasing spatial Cauchy differences over the displayed refinement sequence. 
 
\begin{table}[htbp]
\centering\small 
\setlength{\tabcolsep}{4pt}
\caption{Spatial $L^2$ Cauchy differences at $t=0.02$, $\tau=0.0005$.}
\label{tab:spectral_space}
\begin{tabular}{lrrrr}
\toprule
Variable & $9$ vs. $17$ & $17$ vs. $33$ & $33$ vs. $65$ & $65$ vs. $129$ \\
\midrule
$\varphi$ & $8.602\times10^{-6}$ & $3.119\times10^{-10}$ & $4.105\times10^{-12}$ & $7.828\times10^{-16}$ \\
$u_x$ & $1.131\times10^{-6}$ & $1.088\times10^{-10}$ & $4.679\times10^{-11}$ & $6.474\times10^{-17}$ \\
$u_y$ & $1.763\times10^{-6}$ & $1.298\times10^{-10}$ & $4.676\times10^{-11}$ & $6.232\times10^{-17}$ \\
$\sigma_{xx}$ & $2.752\times10^{-7}$ & $4.180\times10^{-11}$ & $1.141\times10^{-11}$ & $1.890\times10^{-15}$ \\
$\sigma_{xy}$ & $2.463\times10^{-7}$ & $7.740\times10^{-11}$ & $2.294\times10^{-12}$ & $2.370\times10^{-17}$ \\
$\sigma_{yy}$ & $2.846\times10^{-7}$ & $4.189\times10^{-11}$ & $1.141\times10^{-11}$ & $1.960\times10^{-15}$ \\
$\nu$ & $1.178\times10^{-4}$ & $9.336\times10^{-9}$ & $3.368\times10^{-11}$ & $8.564\times10^{-14}$ \\
$P$ & $5.474\times10^{-5}$ & $1.080\times10^{-8}$ & $8.668\times10^{-12}$ & $2.393\times10^{-14}$ \\
\bottomrule
\end{tabular}
\end{table}

\subsection{Unconditional energy-stability verification}

We retain the preceding material parameters and increase the phase,
velocity, and diagonal conformation amplitudes $0.12$, $0.04$, and $0.035$
in \eqref{eq:spectral_initial} to $0.28$, $0.12$, and $0.15$;
the off-diagonal conformation amplitude is $0.075$.
We fix the spatial grid at $129\times129$
with target time $t=1$ and planned step sizes
\begin{equation*}
 \tau=10^{-1},\quad10^{-2},\quad10^{-3},\quad10^{-4}.
\end{equation*}

\begin{figure}[!th]
\centering
\subfigure[Energy evolution.]{\includegraphics[width=0.47\textwidth]{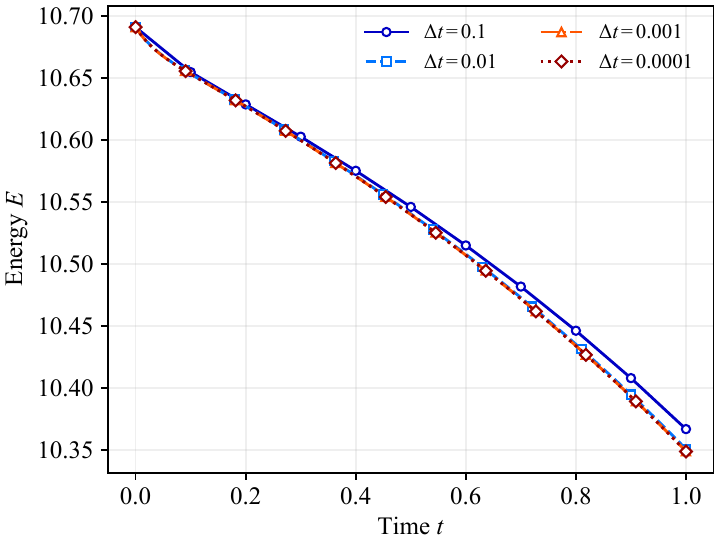}}
\subfigure[Positive definiteness of conformation tensor.]{\includegraphics[width=0.47\textwidth]{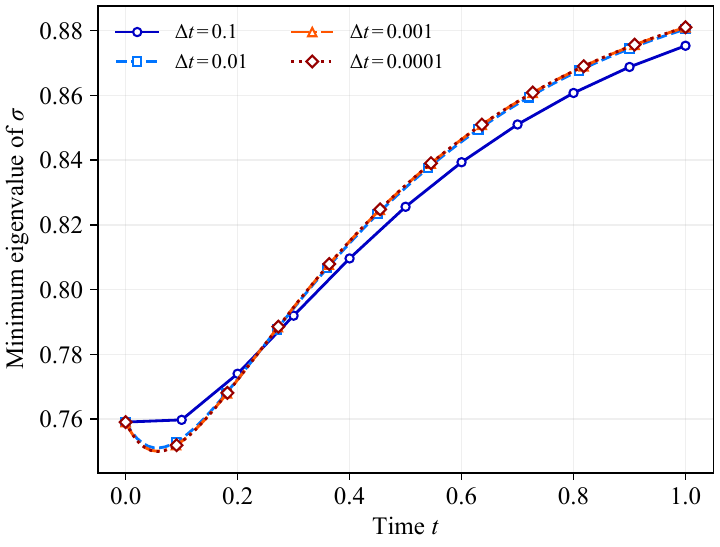}}
\caption{Energy evolution for different time steps with the final time $t=1$.}
\label{fig:energy}
\end{figure}

The energy curves in Figure \ref{fig:energy} share a decreasing trend, 
with differences between the time resolutions. 
The conformation tensor remains positive definite. 
With finer time steps, the smallest conformation eigenvalue first decreases and then recovers.

\subsection{Deformation of droplets in Couette flow}

We denote Newtonian and FENE-P viscoelastic fluids by N and V.
The cases are V in N, N in V, and V in V; $\varphi>0$ is the droplet phase.
All three cases retain identical density, geometry, imposed shear force,
and zero-shear total viscosity. A Newtonian phase has
$\mu_s=1$, $\mu_p=0$ so that $\tau_p$ is zero; a viscoelastic phase has
$\mu_s=0.15$, $\mu_p=0.935$.
Thus $\mu_s+\beta\mu_pc_{\rm eq}=1$ in every bulk phase.
The phase-specific relaxation factors $\lambda^{+}=3$ and $\lambda^{-}=1$
are unchanged between the three cases. Both bulk relaxation factors are positive, since the interpolated factor $\lambda(\varphi)$ occurs in a denominator.

We use $\rho^{+}=1, \rho^{-}=1.5$, $N=129$, $\tau=0.005$, $t=2$,
$\mathrm{Re}=\beta=1$, $\mathrm{Ca}=2$,
$\mathrm{Wi}=2$, $L^{2}=20$, $\delta=0.045$, and $m=0.0005$.
The same nearly circular periodic droplet is initialized in each case:
\begin{equation}
 \varphi_0=\tanh\!\left(\frac{R^2-r_{\rm per}^2}{2\sqrt2\delta R}\right),
 \quad R=0.35,\quad
 r_{\rm per}^2=\frac4{\pi^2}
 \left(\sin^2\!\left(\frac{\pi x}{2}\right)
       +\sin^2\!\left(\frac{\pi y}{2}\right)\right).
\label{profilephase}
\end{equation}

We measure deformation by $D=(a-b)/(a+b)$, where $a\geq b$ are the square roots of the eigenvalues of the droplet's phase-weighted second-moment matrix. The weights use a bounded smooth phase fraction and the coordinate offsets use periodic distances from the droplet centroid. This moment-based measure is used consistently for all pairings.

The domain is periodic in $x$ and
bounded by moving walls at $y=\pm1$, where $\bm u=(P,0)^T$ on $y=1$ and  $\bm u=(-P,0)^T$ on $y=-1$  with $P=1.5$.
The phase field and chemical potential have homogeneous normal-flux
conditions at the walls. 
Initially $\bm u=(Py,0)^T$ and $\bm\sigma=c_{\rm eq}I$.
The initial phase profile uses the same hyperbolic tangent in (\ref{profilephase}) with a different $r_{\rm per}^2$, given by
$r_{\rm per}^2=4\sin^2(\pi x/2)/\pi^2+y^2$. 

\begin{figure}[!th]
\centering
\subfigure[V in N at $t=2$.]{\includegraphics[width=0.32\textwidth]{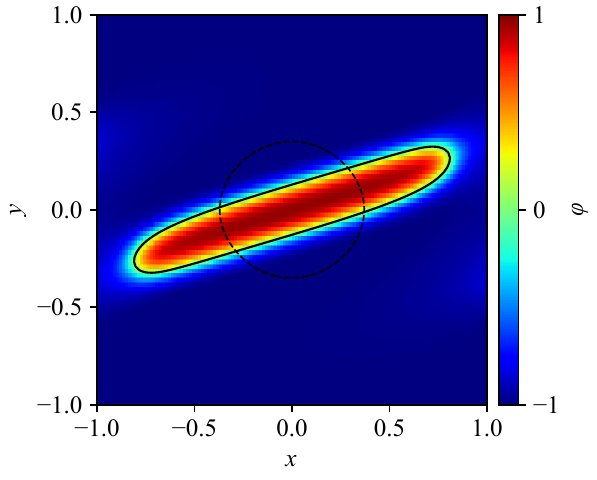}}
\subfigure[N in V at $t=2$.]{\includegraphics[width=0.32\textwidth]{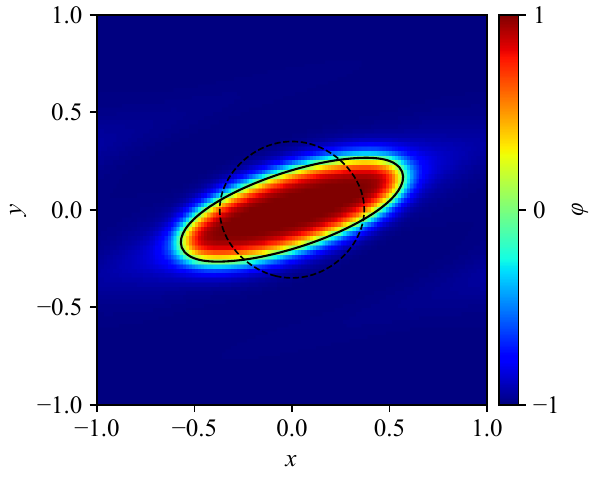}}
\subfigure[V in V at $t=2$.]{\includegraphics[width=0.32\textwidth]{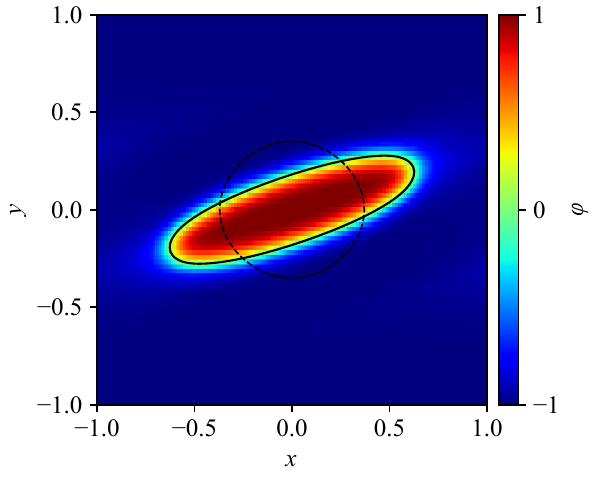}}
\caption{Completed boundary-driven shear cases at $t=2$. Solid and dashed black contours show the final and initial interfaces.}
\label{fig:couette_comparison}
\end{figure}

\begin{figure}[!th]
\centering
\subfigure[Deformation.]{\includegraphics[width=0.47\textwidth]{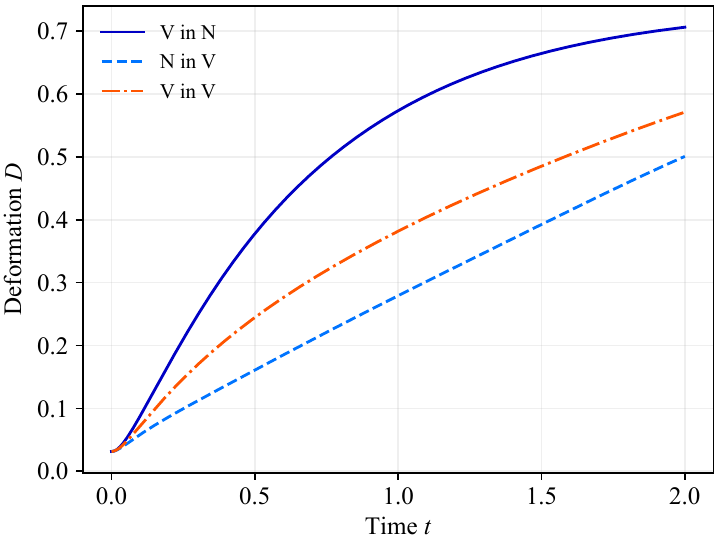}}
\subfigure[Polymer first normal stress difference.]{\includegraphics[width=0.47\textwidth]{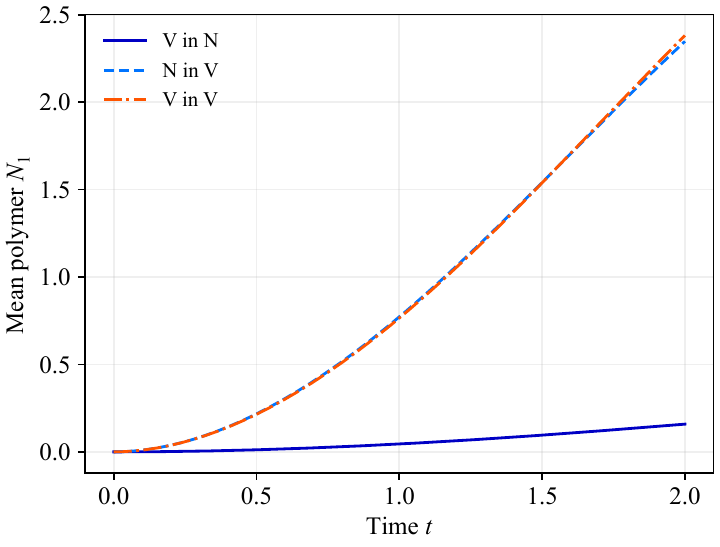}}
\caption{Time histories for the completed cases V in N, N in V, and V in V.}
\label{fig:couette_histories}
\end{figure}

\begin{figure}[!th]
\centering
\subfigure[V in N at $t=0$.]{\includegraphics[width=0.245\textwidth]{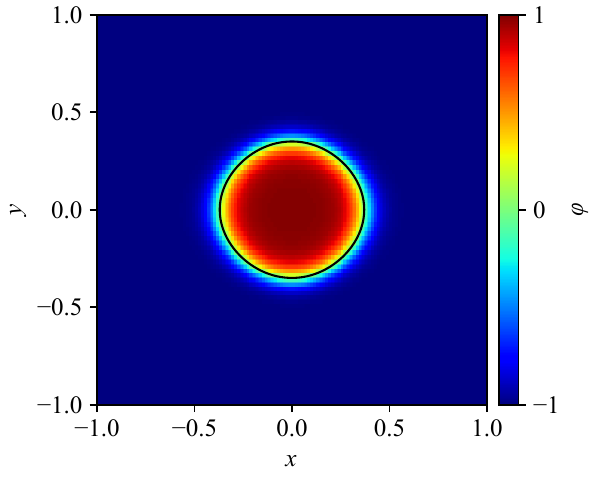}}
\subfigure[V in N at $t=0.5$.]{\includegraphics[width=0.245\textwidth]{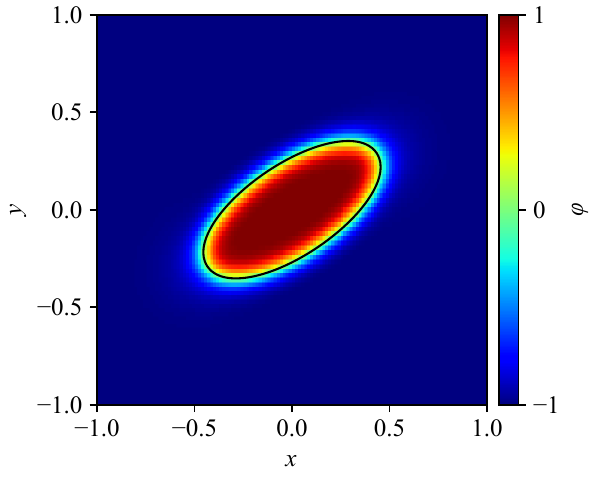}}
\subfigure[V in N at $t=1$.]{\includegraphics[width=0.245\textwidth]{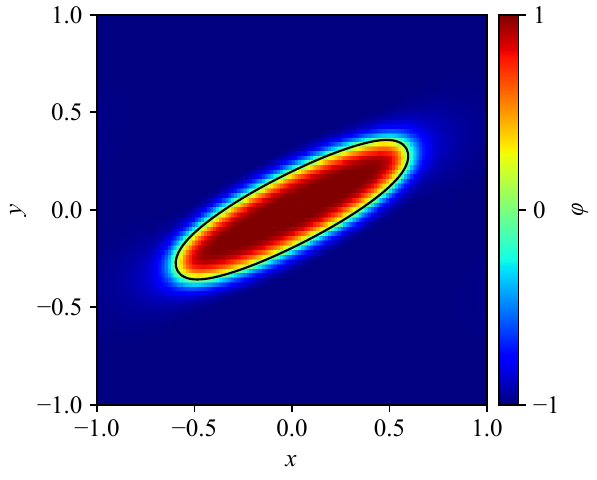}}
\subfigure[V in N at $t=2$.]{\includegraphics[width=0.245\textwidth]{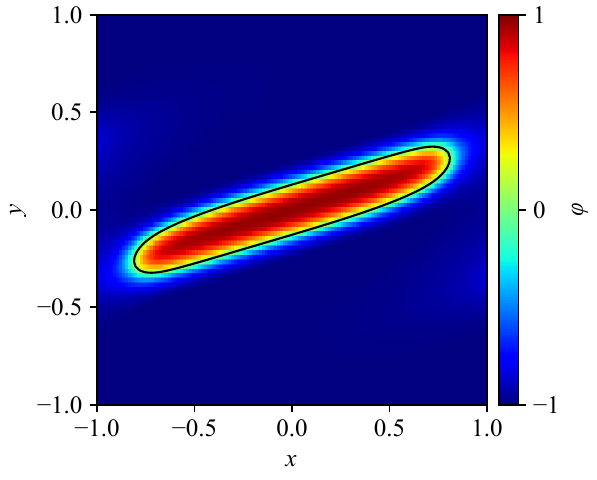}}
\\
\subfigure[N in V at $t=0$.]{\includegraphics[width=0.245\textwidth]{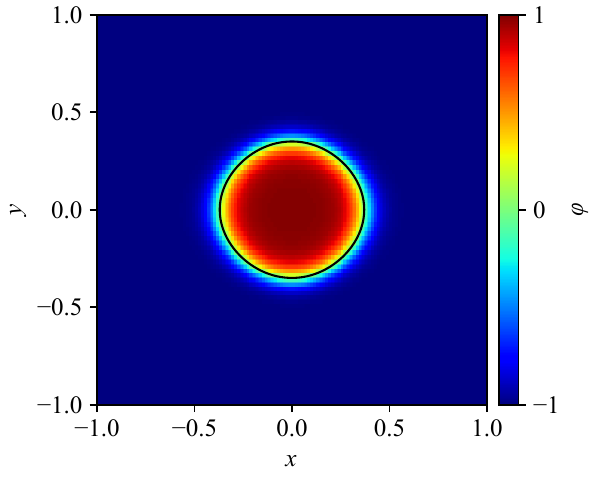}}
\subfigure[N in V at $t=0.5$.]{\includegraphics[width=0.245\textwidth]{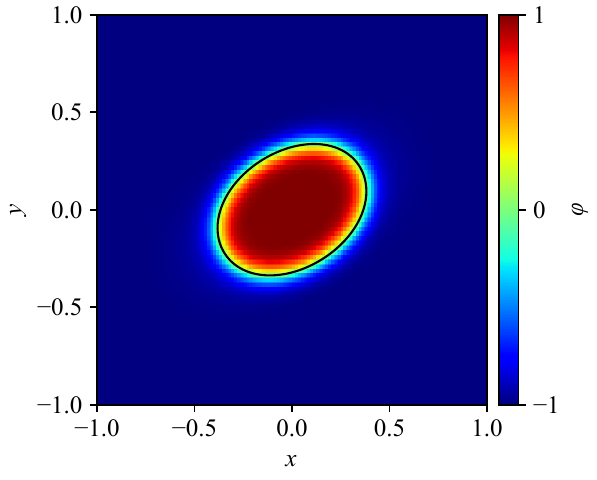}}
\subfigure[N in V at $t=1$.]{\includegraphics[width=0.245\textwidth]{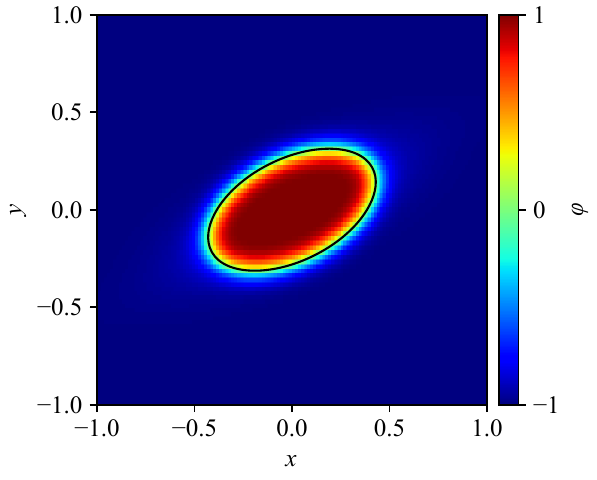}}
\subfigure[N in V at $t=2$.]{\includegraphics[width=0.245\textwidth]{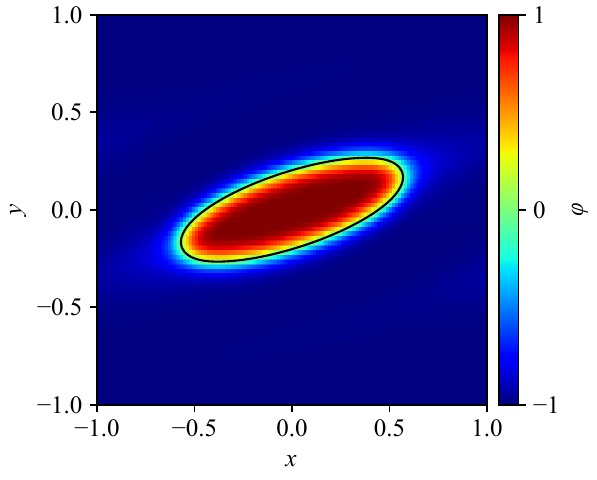}}
\\
\subfigure[V in V at $t=0$.]{\includegraphics[width=0.245\textwidth]{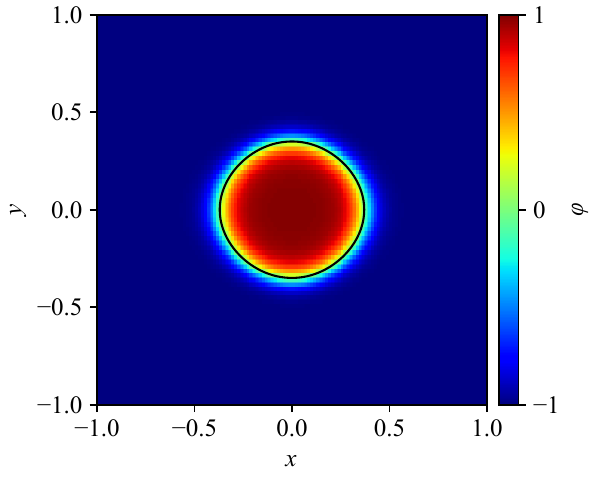}}
\subfigure[V in V at $t=0.5$.]{\includegraphics[width=0.245\textwidth]{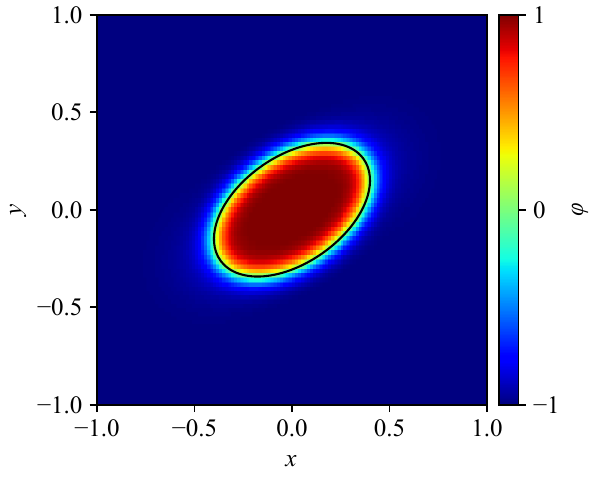}}
\subfigure[V in V at $t=1$.]{\includegraphics[width=0.245\textwidth]{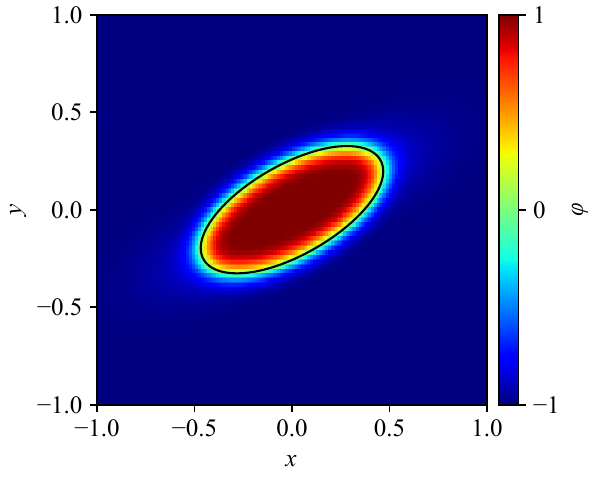}}
\subfigure[V in V at $t=2$.]{\includegraphics[width=0.245\textwidth]{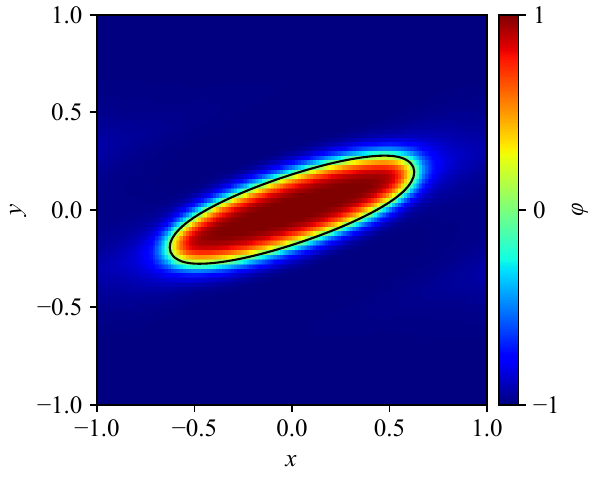}}
\caption{Interface evolution for the completed cases V in N, N in V, and V in V.}
\label{fig:couette_evolution}
\end{figure}

\begin{figure}[!th]
\centering
\subfigure[Energy balance including supplied work.]{\includegraphics[width=0.47\textwidth]{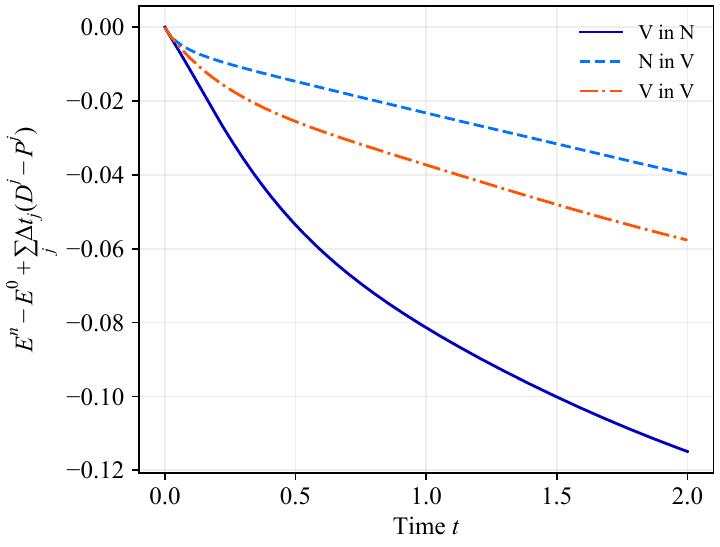}}
\subfigure[Cumulative supplied work.]{\includegraphics[width=0.47\textwidth]{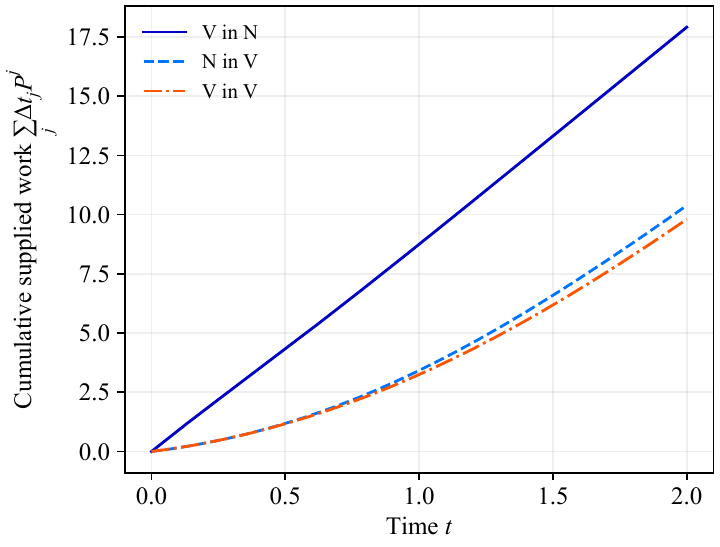}}
\caption{Energy balance and supplied work for the completed cases V in N, N in V, and V in V.}
\label{fig:couette_energy_balance}
\end{figure}

At the final comparison, N in V is the least deformed, V in N is the most deformed, 
and V in V lies between them. The final inclination angles are similar across the fluid pairings.

The Couette snapshots show progressive elongation of a coherent droplet, 
accompanied by increasing deformation and spatially averaged polymer normal stress over the reported interval. 
In the final snapshots, the viscoelastic matrix cases are less deformed than V in N, 
and V in V lies between the two mixed pairings. The ordering differs from the periodic comparison. 
Continued growth of deformation near the end of the recorded interval means 
that the final snapshots should not be interpreted as established steady states.

The plotted normal stress is the domain average of $N_1=({\tau_p})_{xx}-({\tau_p})_{yy}$, 
the polymer contribution to the first normal stress difference in the momentum equation. 
At the final comparison, the ordering of this normal stress across the fluid pairings differs from the ordering of deformation. 
The energy balance includes the supplied moving wall work, so an increase in the driven system's physical energy is compatible with dissipation.

\section{Conclusions}\label{sec7}
We have formulated a thermodynamically consistent two-phase viscoelastic
FENE-P model using the energetic variational approach. A phase field governed
by the Cahn--Hilliard equation describes the diffuse interface and connects
the material properties of the two fluids. The formal matched asymptotic
analysis identifies the corresponding sharp-interface limit as the interface
thickness tends to zero.

For this model, we have constructed a coupled convex-splitting scheme in
time. An exact divided difference accounts for the phase dependence of the
elastic energy, and the time-discrete transport and stress
terms preserve the energy exchange between the phase field, flow, and polymer
conformation. For every admissible, sufficiently regular time-discrete solution, the resulting
energy identity yields an energy estimate without a time-step restriction.
The residual identity quantifies the energy defect introduced by an inexact
nonlinear solve. This analysis establishes energy stability conditional on
admissibility; it does not establish unconditional solvability or convergence
of the nonlinear iteration.

The periodic experiments use Fourier collocation, whereas the
Couette experiments combine Fourier and Legendre collocation.
The smooth test gives temporal Cauchy orders near one and decreasing
spatial Cauchy differences, and the unforced tests show energy decay
for all tested time steps. The wall driven shear experiments compare
viscoelastic droplets and matrices through their deformation and inclination,
with the supplied work included in each energy balance. 

\section*{Acknowledgment}
The research of Y. Qin was partially supported by the National Natural Science Foundation of China (Grant No.12201369). The research of D. He was partially supported by the National Natural Science Foundation of China (Grant No.12672336). 

\section*{Appendix A: Derivation of the two-phase FENE-P model}\label{labelapp2}

Consider the general energy $E$ defined as follows:
\begin{equation}
E = \int_{\Omega} e \mathrm{d} \bm{x}. 
\end{equation}

By using the periodic boundary conditions or the boundary conditions (\ref{BC_zero}), the following identities can be obtained: 
\begin{equation}
\frac{\mathrm{d}E}{\mathrm{d}t} 
= \int_{\Omega} \frac{\partial e}{\partial t} \mathrm{d} \bm{x} 
= \int_{\Omega} \left( \frac{\partial e}{\partial t} 
+ \nabla \cdot \left( \bm{u} e \right) \right) \mathrm{d} \bm{x} 
= \int_{\Omega} \left( \frac{\partial e}{\partial t} 
+ \bm{u} \cdot \nabla e 
+ e \nabla \cdot \bm{u} \right) \mathrm{d} \bm{x} 
= \int_{\Omega} \frac{\mathrm{D} e}{\mathrm{D} t} 
\mathrm{d} \bm{x},
\end{equation} 
where the continuity equation (\ref{eqn:DimNablaU}) is used. 

First, we get the time derivative of kinetic energy as follows. 
\begin{equation}
\begin{aligned}
\frac{\mathrm{d}E_{kin}}{\mathrm{d}t}  
= & \int_{\Omega} \left( \frac{\mathrm{D}\left( \rho \bm{u} \right)}{\mathrm{D}t} \cdot \bm{u} 
- \frac{1}{2} \frac{\mathrm{D}\rho}{\mathrm{D}t} \left| \bm{u} \right|^{2} \right) \mathrm{d}\bm{x} = \int_{\Omega} \left( \left( - \nabla p 
+ \nabla \cdot \bm{\tau}_{s} 
+ \nabla \cdot \bm{\tau}_{p} 
+ \nabla \cdot \bm{\tau}_{\varphi} \right) \cdot \bm{u} 
- \frac{1}{2}\frac{\mathrm{D}\rho}{\mathrm{D}t} 
\left| \bm{u} \right|^{2} \right) \mathrm{d}\bm{x} 
\\ 
= & \int_{\Omega} \left( p \nabla \cdot \bm{u}
- \bm{\tau}_{s} : \nabla \bm{u} 
- \bm{\tau}_{p} : \nabla \bm{u}
+ \left( \nabla \cdot \bm{\tau}_{\varphi} \right) \cdot \bm{u} \right) \mathrm{d}\bm{x} 
- \int_{\Omega} \frac{\rho^{\prime} \left( \varphi \right)}{2} 
\frac{\mathrm{D}\varphi}{\mathrm{D}t} \left| \bm{u} \right|^{2} \mathrm{d}\bm{x} 
\\ 
= & \int_{\Omega} \left( - \mu_{s} \left( \varphi \right) \left( \nabla \bm{u} + \nabla \bm{u}^{T}\right) : \nabla \bm{u} 
- \bm{\tau}_{p} : \nabla \bm{u}
+ \left( \nabla \cdot \bm{\tau}_{\varphi} \right) \cdot \bm{u} \right) \mathrm{d}\bm{x} 
+ \int_{\Omega} \frac{\rho^{\prime} \left( \varphi \right)}{2} 
\left( \nabla \cdot \bm{J} \right) \left| \bm{u} \right|^{2} \mathrm{d}\bm{x} 
\\ 
= & - \int_{\Omega} \frac{\mu_{s} \left( \varphi \right)}{2} \left| \nabla \bm{u} + \nabla \bm{u}^{T}\right|^{2} \mathrm{d}\bm{x} 
+ \int_{\Omega} \left( - \bm{\tau}_{p} : \nabla \bm{u}
+ \left( \nabla \cdot \bm{\tau}_{\varphi} \right) \cdot \bm{u} \right) \mathrm{d}\bm{x} 
+ \int_{\Omega} \frac{\rho^{\prime} \left( \varphi \right)}{2} 
\left| \bm{u} \right|^{2} \nabla \cdot \bm{J} \mathrm{d}\bm{x}. 
\end{aligned}
\end{equation}
Then, we take the derivative of elastic energy shown in the following, \begin{equation}
\begin{aligned}
\frac{\mathrm{d}E_{ela}}{\mathrm{d}t}  
= & - \int_{\Omega}  \frac{\mathrm{D}}{\mathrm{D}t}\left(\frac{\mu_{p}\left( \varphi \right)}{2\lambda\left( \varphi \right)} \right) 
\left( L^{2} \ln \left( 1-\frac{\mathrm{tr}\left(\bm{\sigma}\right)}{L^{2}}\right) 
+ \mathrm{tr} \left(\ln \left( \bm{\sigma}\right) +\bm{I}\right)  \right) 
\mathrm{d} \bm{x}\\
&-\int_{\Omega} 
\frac{\mu_{p}\left( \varphi \right)}{2\lambda\left( \varphi \right)} 
\frac{\mathrm{D}}{\mathrm{D}t} 
\left( L^{2} \ln \left( 1-\frac{\mathrm{tr}\left(\bm{\sigma}\right)}{L^{2}}\right) 
+ \mathrm{tr} \left(\ln \left( \bm{\sigma}\right) +\bm{I}\right)  \right) \mathrm{d} \bm{x} 
\\ 
= & - \int_{\Omega} \left( L^{2} \ln \left( 1-\frac{\mathrm{tr}\left(\bm{\sigma}\right)}{L^{2}}\right) 
+ \mathrm{tr} \left(\ln \left( \bm{\sigma}\right) +\bm{I}\right)  \right) 
\left( \frac{\mu_{p}\left( \varphi \right)}{2\lambda\left( \varphi \right)} \right)^{\prime}
\frac{\mathrm{D}\varphi}{\mathrm{D}t} 
\mathrm{d} \bm{x} 
+ \int_{\Omega} 
\frac{\mu_{p}\left( \varphi \right)}{2\lambda\left( \varphi \right)} 
\frac{\mathrm{D}\bm{\sigma}}{\mathrm{D}t} : 
\left(f \left( \bm{\sigma} \right) \bm{\sigma}^{-1} \right) 
\mathrm{d} \bm{x} 
\\ 
= & \int_{\Omega} \left( \frac{\mu_{p}\left( \varphi \right)}{2\lambda\left( \varphi \right)} \right)^{\prime}
\left( L^{2} \ln \left(1-\frac{\mathrm{tr}\left(\bm{\sigma}\right)}{L^{2}}\right) 
+ \mathrm{tr} \left(\ln \left( \bm{\sigma}\right) +\bm{I}\right) \right) 
\nabla \cdot \bm{J} 
\mathrm{d} \bm{x} 
+ \int_{\Omega} 
\frac{\mu_{p}\left( \varphi \right)}{2\lambda\left( \varphi \right)} 
\frac{\mathrm{D}\bm{\sigma}}{\mathrm{D}t} : 
\left(f \left( \bm{\sigma} \right) \bm{\sigma}^{-1} \right) \mathrm{d} \bm{x}\\
=&\int_{\Omega} \left( \frac{\mu_{p}\left( \varphi \right)}{2\lambda\left( \varphi \right)} \right)^{\prime}
\left( L^{2} \ln \left(1-\frac{\mathrm{tr}\left(\bm{\sigma}\right)}{L^{2}}\right) 
+ \mathrm{tr} \left(\ln \left( \bm{\sigma}\right) +\bm{I}\right) \right) 
\nabla \cdot \bm{J} 
\mathrm{d} \bm{x} \\
&+  \int_{\Omega} 
\frac{\mu_{p}\left( \varphi \right)}{2\lambda\left( \varphi \right)} 
\left( \nabla \bm{u} \cdot \bm{\sigma}
+ \bm{\sigma} \cdot \nabla \bm{u}^{T}
-\frac{1}{\lambda\left(\varphi\right)}f \left( \bm{\sigma} \right) \right) : 
\left(f \left( \bm{\sigma} \right) \bm{\sigma}^{-1} \right) 
\mathrm{d} \bm{x} \\
=&\int_{\Omega} \left( \frac{\mu_{p}\left( \varphi \right)}{2\lambda\left( \varphi \right)} \right)^{\prime}
\left( L^{2} \ln \left(1-\frac{\mathrm{tr}\left(\bm{\sigma}\right)}{L^{2}}\right) 
+ \mathrm{tr} \left(\ln \left( \bm{\sigma}\right) +\bm{I}\right) \right) 
\nabla \cdot \bm{J} 
\mathrm{d} \bm{x} \\
&+ \int_{\Omega} 
\frac{\mu_{p}\left( \varphi \right)}{2\lambda\left( \varphi \right)} 
\frac{2\bm{\sigma}}{1-\frac{\mathrm{tr}\left( \bm{\sigma} \right)}{L^{2}}} : \nabla \bm{u} \mathrm{d}\bm{x} 
- \int_{\Omega} \frac{\mu_{p}\left( \varphi \right)}{2\lambda\left( \varphi \right)} \frac{1}{\lambda\left(\varphi\right)}f \left( \bm{\sigma} \right) : 
\left(f \left( \bm{\sigma} \right) \bm{\sigma}^{-1} \right) 
\mathrm{d} \bm{x} \\
=&\int_{\Omega} \left( \frac{\mu_{p}\left( \varphi \right)}{2\lambda\left( \varphi \right)} \right)^{\prime}
\left( L^{2} \ln \left(1-\frac{\mathrm{tr}\left(\bm{\sigma}\right)}{L^{2}}\right) 
+ \mathrm{tr} \left(\ln \left( \bm{\sigma}\right) +\bm{I}\right) \right) 
\nabla \cdot \bm{J} 
\mathrm{d} \bm{x} \\
&+ \int_{\Omega} 
\frac{\mu_{p}\left( \varphi \right)}{\lambda\left( \varphi \right)} 
\left(\frac{\bm{\sigma}}{1-\frac{\mathrm{tr}\left( \bm{\sigma} \right)}{L^{2}}}-\bm I\right) : \nabla \bm{u} \mathrm{d}\bm{x} 
- \int_{\Omega} \frac{\mu_{p}\left( \varphi \right)}{2\lambda^{2}\left( \varphi \right)} 
\mathrm{tr}\left(f^{2} \left( \bm{\sigma} \right) \bm{\sigma}^{-1} \right) \mathrm{d}\bm{x} 
\end{aligned}
\end{equation} 
The following identity uses the convention
$(\nabla\bm u)_{ij}=\partial_j u_i$ stated in Section~\ref{sec2}.
The logarithmic identity is also used~\cite{Barrett2018FENEP}: 
\begin{equation}\label{eqn:sigma_u}
\left( \nabla \bm{u} \cdot \bm{\sigma}
+ \bm{\sigma} \cdot \nabla \bm{u}^{T} \right) :
\left(f \left( \bm{\sigma} \right) \bm{\sigma}^{-1} \right) 
= 2f(\bm\sigma):\nabla\bm u,\quad \frac{\mathrm{D}  ( \mathrm{tr}(\ln \bm{\sigma}))}{\mathrm{D} t}=\frac{\mathrm{D}  \bm{\sigma}}{\mathrm{D} t}:\bm{\sigma}^{-1}.
\end{equation}
Since $\nabla\cdot\bm u=0$, one also has
$2f(\bm\sigma):\nabla\bm u=2g(\bm\sigma)\bm\sigma:\nabla\bm u$,
which justifies the preceding elastic work simplification.

Finally, the time derivative of mixing energy for the diffuse interface is presented as follows, 
\begin{equation}
\begin{aligned}
\frac{\mathrm{d}E_{mix}}{\mathrm{d}t}  
= & \int_{\Omega} \gamma \frac{\mathrm{D}}{\mathrm{D}t} 
\left( \frac{\delta}{2} \left| \nabla \varphi \right|^{2} 
+ \frac{1}{\delta} W \left( \varphi \right) \right) \mathrm{d}\bm{x} 
\\ 
= & \int_{\Omega} \gamma \left( \delta \nabla \varphi \cdot \nabla \frac{\mathrm{D}\varphi}{\mathrm{D}t}  
- \delta \nabla \varphi \cdot \nabla \bm{u} \cdot \nabla \varphi 
+ \frac{1}{\delta} W^{\prime} \left( \varphi \right) \frac{\mathrm{D}\varphi}{\mathrm{D}t} \right) \mathrm{d}\bm{x} 
\\ 
= & \int_{\Omega} \left( \gamma \left( - \delta \Delta \varphi 
+ \frac{1}{\delta} W^{\prime} \left( \varphi \right) \right) \frac{\mathrm{D}\varphi}{\mathrm{D}t}  
+ \bm{u} \cdot \nabla \cdot \left( \gamma \delta \nabla \varphi \otimes \nabla \varphi \right) \right) \mathrm{d}\bm{x} 
\\ 
= & \int_{\Omega} \left( - \frac{\delta E_{mix}}{\delta \varphi} 
\nabla \cdot \bm{J} + \bm{u} \cdot \nabla \cdot \gamma \delta 
\left( \nabla \varphi \otimes \nabla \varphi \right) \right) 
\mathrm{d}\bm{x}, 
\end{aligned}
\end{equation}
where 
\begin{align*}
\frac{\delta E_{mix}}{\delta \varphi} =\gamma \left( - \delta \Delta \varphi 
+ \frac{1}{\delta} W^{\prime} \left( \varphi \right) \right) .
\end{align*}

Combining the above four equations, we can obtain the following time derivative of the total energy of the FENE-P-CH system 
\begin{equation}
\begin{aligned}
\frac{\mathrm{d}E_{total}}{\mathrm{d}t} 
= & \frac{\mathrm{d}E_{kin}}{\mathrm{d}t} 
+ \frac{\mathrm{d}E_{ela}}{\mathrm{d}t} 
+ \frac{\mathrm{d}E_{mix}}{\mathrm{d}t} 
\\ 
= & - \int_{\Omega} \frac{\mu_{s} \left( \varphi \right)}{2} 
\left| \nabla \bm{u} + \nabla \bm{u}^{T}\right|^{2} \mathrm{d}\bm{x} 
- \int_{\Omega} \bm{\tau}_{p} : \nabla \bm{u} \mathrm{d}\bm{x}
+ \int_{\Omega} \left( \nabla \cdot \bm{\tau}_{\varphi} \right) \cdot \bm{u} \mathrm{d}\bm{x} 
+ \int_{\Omega} \frac{\rho^{\prime} \left( \varphi \right)}{2} 
\left| \bm{u} \right|^{2} \nabla \cdot \bm{J} \mathrm{d}\bm{x} 
\\ 
& + \int_{\Omega} \left( \frac{\mu_{p}\left( \varphi \right)}{2\lambda\left( \varphi \right)} \right)^{\prime}
\left( L^{2} \ln \left(1-\frac{\mathrm{tr}\left(\bm{\sigma}\right)}{L^{2}}\right) 
+ \mathrm{tr} \left(\ln \left( \bm{\sigma}\right) +\bm{I}\right) \right) 
\nabla \cdot \bm{J} 
\mathrm{d} \bm{x} 
\\ 
& + \int_{\Omega} 
\frac{\mu_{p}\left( \varphi \right)}{\lambda\left( \varphi \right)} 
\left(\frac{\bm{\sigma}}{1-\frac{\mathrm{tr}\left( \bm{\sigma} \right)}{L^{2}}}-\bm I\right) : \nabla \bm{u} \mathrm{d}\bm{x} 
- \int_{\Omega} \frac{\mu_{p}\left( \varphi \right)}{2\lambda^{2}\left( \varphi \right)} 
\mathrm{tr}\left(f^{2} \left( \bm{\sigma} \right) \bm{\sigma}^{-1} \right) \mathrm{d}\bm{x} \\
&+ \int_{\Omega} \left( - \frac{\delta E_{mix}}{\delta \varphi} \nabla \cdot \bm{J} 
+ \bm{u} \cdot \nabla \cdot \left( \gamma \delta \nabla \varphi \otimes \nabla \varphi \right) \right) \mathrm{d}\bm{x} 
\\ 
= & - \int_{\Omega} \frac{\mu_{s} \left( \varphi \right)}{2} 
\left| \nabla \bm{u} + \nabla \bm{u}^{T}\right|^{2} \mathrm{d}\bm{x} 
+ \int_{\Omega} \bm{u} \cdot \nabla \cdot \left(\bm{\tau}_{\varphi} 
+ \gamma \delta \nabla \varphi \otimes \nabla \varphi 
+ \rho^{\prime} \left( \varphi \right) \bm{u} \otimes \bm{J} \right) \mathrm{d}\bm{x} 
\\ 
& + \int_{\Omega} \left( \left( \frac{\mu_{p}\left( \varphi \right)}{2\lambda\left( \varphi \right)} \right)^{\prime}
\left( L^{2} \ln \left(1-\frac{\mathrm{tr}\left(\bm{\sigma}\right)}{L^{2}}\right) 
+ \mathrm{tr} \left(\ln \left( \bm{\sigma}\right) +\bm{I}\right) \right) 
- \frac{\delta E_{mix}}{\delta \varphi} \right) \nabla \cdot \bm{J} 
\mathrm{d} \bm{x} 
\\ 
& + \int_{\Omega} 
\left( 
\frac{\mu_{p}\left( \varphi \right)}{\lambda\left( \varphi \right)} \left(\frac{\bm{\sigma}}{1-\frac{\mathrm{tr}\left( \bm{\sigma} \right)}{L^{2}}}-\bm I\right) -\bm{\tau}_{p} \right) : \nabla \bm{u} \mathrm{d}\bm{x} 
- \int_{\Omega} \frac{\mu_{p}\left( \varphi \right)}{2\lambda^{2}\left( \varphi \right)} 
\mathrm{tr}\left(f^{2} \left( \bm{\sigma} \right) \bm{\sigma}^{-1} \right) 
\mathrm{d} \bm{x}, 
\end{aligned}
\end{equation}
where we use the fact that the density $\rho$ 
is a linear combination of the densities in the two phases, 
$\rho^{\prime}\left(\varphi\right)$ is a constant, 
which means 
\begin{equation}
\int_{\Omega} \frac{\rho^{\prime} \left( \varphi \right)}{2} 
\left| \bm{u} \right|^{2} \nabla \cdot \bm{J} \mathrm{d} \bm{x} 
= - \int_{\Omega} \rho^{\prime} \left( \varphi \right) 
\bm{u} \cdot \left( \nabla \bm{u} \cdot \bm{J} \right) \mathrm{d} \bm{x} 
= \int_{\Omega} \bm{u} \cdot \left( \nabla \cdot \left( \rho^{\prime} \left( \varphi \right) \bm{u} \otimes \bm{J} \right) \right) \mathrm{d} \bm{x}. 
\end{equation}
Therefore we can define 
\begin{subequations}
\begin{align}
\bm{J} = & M\left( \varphi \right) \nabla 
\left( \left( \frac{\mu_{p}\left( \varphi \right)}{2\lambda\left( \varphi \right)} \right)^{\prime}
\left( L^{2} \ln \left(1-\frac{\mathrm{tr}\left(\bm{\sigma}\right)}{L^{2}}\right) 
+ \mathrm{tr} \left(\ln \left( \bm{\sigma}\right) +\bm{I}\right) \right) 
- \frac{\delta E_{mix}}{\delta \varphi} \right)=-M(\varphi)\nabla \nu, 
\\ 
\bm{\tau}_{\varphi} = & - \gamma \delta \nabla \varphi \otimes \nabla \varphi 
- \rho^{\prime} \left( \varphi \right) \bm{u} \otimes \bm{J}
= - \gamma \delta \nabla \varphi \otimes \nabla \varphi 
+ \rho^{\prime} \left( \varphi \right) \bm{u} \otimes 
M(\varphi)\nabla \nu, 
\\ 
\bm{\tau}_{p} = & \frac{\mu_{p}\left( \varphi \right)}{\lambda\left( \varphi \right)} 
f(\bm {\sigma}), 
\end{align}
\end{subequations} 
where
\begin{equation}
\nu=\frac{\delta E_{mix}}{\delta \varphi}-\left( \frac{\mu_{p}\left( \varphi \right)}{2\lambda\left( \varphi \right)} \right)^{\prime}
\left( L^{2} \ln \left(1-\frac{\mathrm{tr}\left(\bm{\sigma}\right)}{L^{2}}\right) 
+ \mathrm{tr} \left(\ln \left( \bm{\sigma}\right) +\bm{I}\right) \right) 
\end{equation}
is the chemical potential. 

\bibliographystyle{unsrt}
\bibliography{refs}
\end{document}